# HIGHER RANK CASE OF DWORK'S CONJECTURE

DAQING WAN

*Dedicated to the memory of Bernard Dwork*

## 1. Introduction

In this series of two papers, we prove the $p$-adic meromorphic continuation of the pure slope L-functions arising from the slope decomposition of an overconvergent F-crystal, as conjectured by Dwork [6]. More precisely, we prove a suitable extension of Dwork's conjecture in our more general setting of $\sigma$-modules, see section 2 for precise definitions of the various notions used in this introduction. Our main result is the following theorem.

**Theorem 1.1.** *Let $X$ be a smooth affine variety defined over a finite field $\mathbf{F}_q$ of characteristic $p > 0$. Let $(M, \phi)$ be a finite rank overconvergent $\sigma$-module over $X/\mathbf{F}_q$. Then, for each rational number $s$, the pure slope $s$ L-function $L_s(\phi, T)$ attached to $\phi$ is $p$-adic meromorphic everywhere.*

The proof of this theorem will be completed in two papers. In the present higher rank paper, we introduce a reduction approach which reduces Theorem 1.1 to the special case when the slope $s$ $(s = 0)$ part of $\phi$ has rank one and the base space $X$ is the simplest affine space $\mathbf{A}^n$. This part is essentially algebraic. It depends on Monsky's trace formula, Grothendieck's specialization theorem, the Hodge-Newton decomposition and Katz's isogeny theorem. In our next paper [23], we will handle the rank one case over the affine space $\mathbf{A}^n$. The rank one case is very much analytic in nature and forces us to work in a more difficult infinite rank setting, generalizing and improving the limiting approach introduced in [19].

Dwork's conjecture grew out of his attempt to understand the $p$-adic analytic variation of the pure pieces of the zeta function of a variety when the variety moves through an algebraic family. To give an important geometric example, let us consider the case that $f: Y \rightarrow X$ is a smooth and proper morphism over $\mathbf{F}_q$ with a smooth and proper lifting to characteristic zero. Berthelot's result [1] says that the relative crystalline cohomology $R^i f_{\mathrm{crys},*} \mathbf{Z}_p$ modulo torsion is an overconvergent F-crystal $M_i$ over $X$. Applying Theorem 1.1, we conclude that the pure L-functions arising from these geometric overconvergent F-crystals $M_i$ are $p$-adic meromorphic. In particular, this implies the existence of an exact $p$-adic formula for geometric $p$-adic character sums and a suitable $p$-adic equi-distribution theorem for the roots of zeta functions. For more detailed arithmetic motivations and further open problems, see the expository papers [20][21].

1991 *Mathematics Subject Classification.* Primary 11G40, 11S40; Secondary 11M41, 14G15.
This work was partially supported by NSF.

An interesting archimedian analogue in the $\ell$-adic ($\ell \neq p$) setting of Dwork's conjecture is closely related to Deligne's conjecture on the purity structure of $\ell$-adic sheaves. Instead of using the slope (or $p$-adic absolute value) in the $p$-adic case, one uses the weight (or complex absolute value) in the $\ell$-adic case to decompose the Euler factors fibre by fibre. Given a constructible $\ell$-adic sheaf $\mathcal{E}$ on $X$, one can similarly attach [21] an L-function $L_s(\mathcal{E}, T)$ to the pure weight $s$ part of $\mathcal{E}$ for each real number $s$. The archimedian version of Dwork's conjecture in the $\ell$-adic setting then says that the pure weight L-function $L_s(\mathcal{E}, T)$ is rational for every $s$. In comparison, Deligne's conjecture says that every constructible $\ell$-adic sheaf is mixed in the sense of [3], i.e., a successive extension of pure sheaves. This conjecture of Deligne together with Grothendieck's rationality theorem immediately implies the above archimedian Dwork conjecture. In the case that the $\ell$-adic sheaf $\mathcal{E}$ is geometric, i.e., arising from the $\ell$-adic higher direct image with compact support, Deligne's main theorem [3] says that such a geometric sheaf $\mathcal{E}$ is indeed mixed. This together with Grothendieck's rationality theorem shows that the archimedian Dwork conjecture is indeed true in the geometric case, but still open in general and depends on Deligne's conjecture.

We now give a brief description of the content of this paper and the main steps of our proof. Section 2 sets up the basic definitions of a finite rank $\sigma$-module $(M, \phi)$ and its L-function $L(\phi, T)$ on a smooth connected affine variety $X$ over $\mathbf{F}_q$, with an emphasis on overconvergent $\sigma$-modules. Restricting each Euler factor of the L-function to its slope $s$ part leads to the pure slope L-function $L_s(\phi, T)$. Dwork's conjecture in our setting then says that the pure slope L-function $L_s(\phi, T)$ is meromorphic if the initial $\sigma$-module $(M, \phi)$ is overconvergent.

Section 3 describes the generalized Monsky trace formula for the L-function of an overconvergent $\sigma$-module. This trace formula was proved by Monsky [13] in the case that $\phi$ has rank one. For the sake of completeness, a proof of the trace formula in the general case is included in the appendix (section 10) of this paper, based on various results in [13]. This trace formula expresses the L-function $L(\phi, T)$ as a finite alterating product of the Fredholm determinants of several completely continuous operators. More precisely, we have

**Theorem 1.2.** *Let $\phi$ be an overconvergent finite rank $\sigma$-module over $X/\mathbf{F}_q$. Then, there are $n+1$ completely continuous operators $\Theta_i(\phi)$ ($0 \leq i \leq n$) acting on certain $p$-adic Banach spaces such that*

$$L(\phi, T) = \prod_{i=0}^{n} \det(I - T\Theta_i(\phi))^{(-1)^{i-1}}.$$

By Serre's $p$-adic spectral theory [16], the Fredholm determinant of a completely continuous operator is entire. Thus, from Theorem 1.2, we deduce the following result.

**Theorem 1.3.** *Let $\phi$ be a finite rank overconvergent $\sigma$-module over $X/\mathbf{F}_q$. Then, the L-function $L(\phi, T)$ is meromorphic.*

Section 4 introduces the notions of basis filtration and basis polygon attached to a $\sigma$-module with respect to a given basis. This makes it possible to define ordinary $\sigma$-modules and study the Hodge-Newton decomposition for an ordinary $\sigma$-module. The Hodge-Newton decomposition identifies the pure slope L-function $L_s(\phi, T)$ in the ordinary case with the L-function $L(\phi_s, T)$ of a new pure slope $\sigma$-module $\phi_s$, where $\phi_s$ is the pure slope $s$ part of $\phi$. That is,

$$L_s(\phi, T) = L(\phi_s, T). \tag{1.1}$$

If the new $\sigma$-module $\phi_s$ were overconvergent, then we would be done by Theorem 1.3. The difficulty is that the new $\sigma$-module $\phi_s$ is no longer overconvergent.

In section 5, using various tensor and exterior constructions, we introduce a reduction approach by proving the meromorphic continuation of $L(\phi_s, T)$ in the easier case when $\phi$ is ordinary and the determinant $\sigma$-module $\det(\phi_s)$ of each pure piece $\phi_s$ of $\phi$ is overconvergent. This part is quite simple and independent of our limiting approach. From formula point of view, our proof amounts to the following result.

**Theorem 1.4.** *Let $\phi$ be a finite rank ordinary overconvergent $\sigma$-module over $X/\mathbf{F}_q$. Assume that for each non-trivial slope $s$, the determinant $\det(\phi_s)$ of the pure slope $s$ part $\phi_s$ of $\phi$ is overconvergent. Then, for each non-trivial slope $s$, there is an infinite sequence $\psi_{j,s}$ $(1 \leq j < \infty)$ of finite rank overconvergent ordinary $\sigma$-modules over $X/\mathbf{F}_q$ such that*

$$L(\phi_s, T) = \prod_{j=1}^{\infty} L(\psi_{j,s}, p^j T)^{\pm 1}. \tag{1.2}$$

Although the product in (1.2) is an infinite product, the key point is the existence of the shifting factor $p^j$ and the overconvergence of $\psi_{j,s}$. The crucial shifting factor shows that the product in (1.2) is really finite in nature if we restrict our attention to the zeros and poles in a finite disk. By Theorem 1.3, we deduce from (1.1) and (1.2) that the pure L-function $L_s(\phi, T)$ is indeed meromorphic if $\phi$ satisfies the conditions of Theorem 1.4.

In section 6, we study the general ordinary case. Twisting by a rank one unit root $\sigma$-module from a suitable Hodge-Newton decomposition, we show that the above reduction approach still works provided one can prove the meromorphic continuation of $L(\phi_s, T)$ in the special case that $\phi_s$ has rank one. Essentially, what we prove here is the following generalization of Theorem 1.4.

**Theorem 1.5.** *Let $\phi$ be a finite rank ordinary overconvergent $\sigma$-module over $X/\mathbf{F}_q$. Then, for each non-trivial slope $s$, there is an infinite sequence $\psi_{j,s}$ $(1 \leq j < \infty)$ of finite rank ordinary $\sigma$-modules over $X/\mathbf{F}_q$ such that*

$$L(\phi_s, T) = \prod_{j=1}^{\infty} L(\psi_{j,s}, p^j T)^{\pm 1} \tag{1.3}$$

*and such that each $\psi_{j,s}$ is of the form*

$$\psi_{j,s} = \Psi_{j,s} \otimes \rho_{j,s}^{-1}, \tag{1.4}$$

*where $\Psi_{j,s}$ is overconvergent, $\rho_{j,s}$ is of rank one and is the slope zero part of some finite rank overconvergent ordinary $\sigma$-module.*

Note that $\rho_{j,s}$ is in general not overconvergent. Thus, we cannot conclude the meromorphic continuation of $L(\phi_s, T)$ yet. By Theorem 1.5, we still need to handle the L-function of the $\sigma$-module in (1.4) which is not overconvergent. The $\sigma$-module in (1.4) is the tensor product of a harmless overconvergent $\sigma$-module and a non-overconvergent but rank one $\sigma$-module. Thus, we have reduced the ordinary case of Dwork's conjecture from higher rank $\phi_s$ to the following slightly more general version of rank one case, called the key lemma.

**Theorem 1.6** (key lemma). *Let $\phi$ be a finite rank overconvergent $\sigma$-module over $X/\mathbf{F}_q$. Let $\rho$ be the slope zero part of some finite rank overconvergent $\sigma$-module ordinary at the slope zero side. Assume that $\rho$ has rank one. Then the L-function $L(\phi \otimes \rho^k, T)$ is meromorphic everywhere for every integer $k$. Furthermore, the family of L-functions $L(\phi \otimes \rho^k, T)$ parametrized by integers $k$ in any given residue class modulo $(q-1)$ is a strong family in the sense of* [23].

In section 9, we further reduce the key lemma from a smooth affine variety $X$ to the affine $n$-space $\mathbf{A}^n$, using some ideas from our limiting approach. The key lemma for the affine space $\mathbf{A}^n$ was proved in [19] for a special type of Frobenius lifting using the infinite rank version of Dwork's trace formula. The full key lemma for the affine $n$-space $\mathbf{A}^n$ will be proved in our next paper [23], using the infinite rank version of the more difficult Monsky trace formula.

In section 7, we turn to studying the harder non-ordinary case in the context of Grothendieck's specialization theorem and Katz's isogeny theorem [10]. Making extensive use of various basic results in [10], we show the above reduction method can be modified to treat the non-ordinary case as well, but we will need a similar non-ordinary key lemma which is a little harder than the ordinary key lemma (Theorem 1.6). We then show that the non-ordinary key lemma can be reduced to the ordinary key lemma, by establishing a crucial stability theorem of the overconvergence under certain weak isogeny map. This then finishes the proof of Theorem 1.1.

In section 8, we treat Dwork's conjecture in full generality by working with his $k$-th power pure slope L-function. This generalization is easily done by using the following decomposition formula (see [19]) expressing the $k$-th power $\phi^k$ (iterate) in terms of symmetric and exterior power products:

$$(1.5) \qquad L(\psi \otimes \phi^k, T) = \prod_{j \geq 1} L(\psi \otimes \mathrm{Sym}^{k-j}\phi \otimes \wedge^j \phi, T)^{(-1)^{j-1}j}.$$

In the easier case of Theorem 1.4, equation (1.2) gives an interesting product formula for the pure slope L-function $L_s(\phi, T)$ in terms of L-functions of finite rank overconvergent $\sigma$-modules. One could ask if there is a similar infinite product formula in the general case of Theorem 1.1. By Theorem 1.5, the key is to obtain such a product formula in the rank one case of the key lemma. The answer is fortunately positive, but the product formula will have to involve L-functions of certain infinite rank overconvergent $\sigma$-modules, called nuclear overconvergent $\sigma$-modules. This is done in [23], which is basically our proof of the key lemma. For the sake of simplicity, in the present paper, we stay in the algebraic finite rank setting.

Eventually, we hope to extend our main results to the infinite rank setting, i.e., when the ambient $\phi$ is an overconvergent nuclear $\sigma$-module of infinite rank, see [23] for the nuclear notion and the infinite rank setup. Currently, our results in this paper and in [23] show that Theorem 1.1 is already true for infinite rank nuclear overconvergent $\phi$ provided that $\phi$ is ordinary. The non-ordinary case requires a substantial amount of additional foundational work.

**Acknowledgement**. The author wishes to thank P. Berthelot, R. Coleman and S.T. Yau for related interesting discussions. The author is also grateful to the referee and the editors (B. Sturmfels and A.J. de Jong) for their attempts in making the paper more readable.

## 2. $\sigma$-MODULES AND L-FUNCTIONS

In this section, we set up the basic definitions that are needed to define the L-function and formulate Dwork's conjecture in a suitable context. Our discussion will be brief, focusing only on what we need from L-function point of view.

Let $R$ be a complete discrete valuation ring of characteristic zero whose residue field $R/\pi R$ is the finite field $\mathbf{F}_q$ of $q$ elements, where $\pi$ is a fixed uniformizer of $R$. Let $n$ be a fixed positive integer and let $\mathbf{N}$ be the set of non-negative integers. Let $R[X_1, \cdots, X_n]^\dagger$ be the overconvergent power series ring consisting of all power series
$$\sum_{u \in \mathbf{N}^n} a_u X^u, a_u \in R$$
satisfying
$$\lim_{k \to \infty} \inf_{|u| > k} \frac{\mathrm{ord}_\pi a_u}{|u|} > 0,$$
where for a lattice point $u = (u_1, \cdots, u_n) \in \mathbf{N}^n$, we define $X^u = X_1^{u_1} \cdots X_n^{u_n}$ and $|u| = u_1 + \cdots + u_n$. This overconvergent ring is the so-called weak completion of the polynomial ring $R[X_1, \cdots, X_n]$ which corresponds to the affine $n$-space $\mathbf{A}^n$ over $R$ or the weak lifting to characteristic zero of the affine $n$-space $\mathbf{A}^n$ over $\mathbf{F}_q$. This ring $R[X_1, \cdots, X_n]^\dagger$ is Noetherian, as shown by Fulton [8], see also [17] for more general results in this direction, generalizing Hilbert's basis theorem.

To treat a more general smooth affine variety (scheme), we let $A$ be a homomorphic image of the commutative $R$-algebra $R[X_1, \cdots, X_n]^\dagger$. Since $R[X_1, \cdots, X_n]^\dagger$ is Noetherian, there are finitely many elements $f_1, \cdots, f_m$ in $R[X_1, \cdots, X_n]^\dagger$ such that
$$A = R[X_1, \cdots, X_n]^\dagger / (f_1, \cdots, f_m).$$
Assume that $A$ is flat over $R$. The ring $A$ is also weakly complete. Its reduction modulo $\pi$ corresponds to an affine variety $X$ defined over $\mathbf{F}_q$. Let $\sigma$ (if it exists) denote a fixed $R$-algebra endomorphism of $A$ such that the reduction of $\sigma$ modulo $\pi$ is the $q$-th power Frobenius map. Namely, $\sigma$ is a lifting to $A$ of the $q$-th power Frobenius map acting on the coordinate ring of the characteristic $p$ variety $X/\mathbf{F}_q$. For an element $a \in A$, the image $\sigma(a)$ of $a$ under the $\sigma$-action will sometimes be denoted by $a^\sigma$, depending on convenience. If the affine variety $X$ over $\mathbf{F}_q$ is smooth, one can always lift the coordinate ring of X to some $A$ (unique up to an isomorphism) flat over $R$ such that the $q$-th power Frobenius map of $X/\mathbf{F}_q$ lifts to an $R$-algebra endomorphism $\sigma$ of $A$. We assume that we are in this situation as in the theory of Monsky-Washnitzer, see [15] for an exposition where the original very smooth condition is improved to the usual smooth condition. Thus, throughout this paper, we shall assume that $A$ is flat over $R$, the reduction $A/\pi A$ is smooth over $\mathbf{F}_q$ and $\sigma$ is a fixed $R$-algebra endomorphism of $A$ lifting the $q$-th power Frobenius of $A/\pi A$. In particular, we always have $\sigma(\pi) = \pi$ since $\sigma$ is an $R$-algebra endomorphism. Without loss of generality, we shall also assume that the smooth affine $\mathbf{F}_q$-variety $X = \mathrm{Spec}(A/\pi A)$ is connected and thus $A$ is an integral domain. The map $\sigma : A \to \sigma(A)$ is clearly bijective. The image $\sigma(A)$ is a subring of $A$. Viewed as a $\sigma(A)$-module, the ring $A$ is a locally free module of rank $q^{\dim(X)}$. Our main concern is the meromorphic continuation of an L-function which is a local property. Thus, shrinking $X$ if necessary, we can assume that $A$ is actually a finite free $\sigma(A)$-module of rank $q^{\dim(X)}$.

Let $A_0$ be the $\pi$-adic completion of $A$, which is also uniquely defined up to an isomorphism. We use the notation $A_0$ instead of $\hat{A}$ to be compatible with our notation $A_c$ of $c$ log-convergent ring, see [20]. Thus, $A_0$ is NOT the reduction of $A$ modulo $\pi$ as some authors would normally mean. The $R$-algebra endomorphism $\sigma$ on $A$ extends easily and uniquely to its completion $A_0$ by continuity. In this paper, we shall only consider $\sigma$-modules of finite rank, although some of our results can be extended directly to infinite rank setting.

**Definition 2.1.** A $\sigma$-module over $A_0$ (or over X) is a pair $(M, \phi)$, where $M$ is a free $A_0$-module of finite rank and $\phi$ is a $\sigma$-linear map:
$$\phi : M \longrightarrow M, \ \phi(ae) = \sigma(a)\phi(e), \ \ a \in A_0, e \in M.$$
An overconvergent $\sigma$-module over $A$ (or over X) is a pair $(M, \phi)$, where $M$ is a free $A$-module of finite rank and $\phi$ is a $\sigma$-linear map:
$$\phi : M \longrightarrow M, \ \phi(ae) = \sigma(a)\phi(e), \ \ a \in A, e \in M.$$

Thus, a $\sigma$-module over $X/\mathbf{F}_q$ is just a finite free $A_0$-module $M$ endowed with a Frobenius action extending the action of $\sigma$ from $A_0$ to $M$. An overconvergent $\sigma$-module over $X/\mathbf{F}_q$ is a finite free $A$-module $M$ endowed with a Frobenius action extending the action of $\sigma$ from $A$ to $M$. For example, the pair $(A_0, \sigma)$ is the identity $\sigma$-module of rank one. The pair $(A, \sigma)$ is the identity overconvergent $\sigma$-module of rank one. After the base extension $A \to A_0$, an overconvergent $\sigma$-module is naturally a $\sigma$ module. We shall often call a $\sigma$-module $(M, \phi)$ to be overconvergent if it can be obtained from an overconvergent $\sigma$-module by base extension $A \to A_0$.

The $\sigma$-linear map $\phi$ in a $\sigma$-module $(M, \phi)$ can also be viewed as an $A_0$-linear map $\phi : M^{(\sigma)} \to M$ of $A_0$-modules, where $M^{(\sigma)} = A_0 \otimes_\sigma M$ is the pull back of $M$ under $\sigma : A_0 \to A_0$. That is, for $a \in A_0$ and $m \in M$, we have the relation
$$\sigma(a) \otimes m = am.$$
The map $\phi$ is called the Frobenius map of the $\sigma$-module $M$. We shall often just write $M$ or $\phi$ for the pair $(M, \phi)$. A $\sigma$-module $(M, \phi)$ is called a unit root $\sigma$-module if the map $\phi : M^{(\sigma)} \to M$ is an isomorphism of $A_0$-modules. A morphism of $\sigma$-modules $f : (M, F) \to (M', F')$ is an $A_0$-linear map such that $fF = F'f$. The category of $\sigma$-modules up to isogeny is obtained from the category of $\sigma$-modules by keeping the same objects, but tensoring the Hom groups with the quotient field of $R$. An isogeny between $\sigma$-modules is a morphism of $\sigma$-modules which becomes an isomorphism in this new category. All these notions carry over to the subcategory of overconvergent $\sigma$-modules with $A_0$ replaced by $A$.

Note that we do not assume that $\phi$ is injective. Nor do we assume that $\text{Coker}(\phi)$ is killed by a power of $\pi$ as in the theory of F-crystals. Thus, it may happen that $\det\phi$ is the zero map while $\phi$ is non-zero. One can define various tensor, symmetric and exterior operations of $\sigma$-modules (resp. overconvergent $\sigma$-modules) in usual way. By convention, we define $\wedge^0 \phi$ to be the rank one unit root $\sigma$-module $(A_0, \sigma)$. It is the identity element in the tensor category of $\sigma$-modules. Similarly, if $\phi$ is overconvergent, we define $\wedge^0 \phi$ to be $(A, \sigma)$ which is the identity element in the tensor category of overconvergent $\sigma$-modules.

*Remark* 2.2. In Definition 2.1, one could use the slightly more general notion of projective $A_0$-modules (or even just finite $A_0$-modules) instead of free $A_0$-modules. This generalization is not essential for our L-function purpose because of a standard reduction procedure. Such a generalization is however, convenient in studying certain properties. Thus, we shall allow $M$ to be a projective module or even a finite $A_0$-module in some explicitly stated cases. This occurs in section 3 and in the appendix. In all other sections, our $M$ is simply free. Our main concern of the present paper is Dwork's conjecture which is local in nature. Thus, for our application, no generality at all is lost by restricting to free $A_0$-modules.

For a given $\sigma$-module $(M, \phi)$ over $X/\mathbf{F}_q$, one can attach an L-function to it. For this purpose, we need to recall the notion of the Teichmüller lifting of a geometric point $x \in X(\bar{\mathbf{F}}_q)$ defined with respect to $\sigma$. Let $\bar{A}$ be the reduction of $A_0$ modulo $\pi$, which is the coordinate ring of $X$ over $\mathbf{F}_q$. For a positive integer $d$, let $R_d$ denote the unramified extension of $R$ with residue field $\mathbf{F}_{q^d}$. Let $\bar{x} \in X$ be a geometric point of degree $d$ over $\mathbf{F}_q$. That is, the point $\bar{x}$ is a surjective $\mathbf{F}_q$-algebra homomorphism in $\operatorname{Hom}_{\mathbf{F}_q}(\bar{A}, \mathbf{F}_{q^d})$. By a theorem of Monsky-Tate, there is a unique surjective $R$-algebra homomorphism $x \in \operatorname{Hom}_R(A, R_d)$ whose reduction mod $\pi$ is $\bar{x}$ such that

$$x \circ \sigma^d = x.$$

The element $x$ is called the Teichmüller lifting of $\bar{x}$ defined with respect to $\sigma$. Let $\tau$ be the Frobenius automorphism of $R_d$ which fixes $\pi$ such that the reduction map $\tau$ modulo $\pi$ is the $q$-th power map of $\mathbf{F}_{q^d}$. One checks that both $x \circ \sigma$ and $\tau \circ x$ are the Teichmüller lifting of $\bar{x}^q$. It follows that

$$x \circ \sigma = \tau \circ x.$$

Thus, the action of $\sigma$ on the Teichmüller lifting $x$ is compatible with the action of Frobenius $\tau$ on $R_d$. To save notation, we shall denote $\tau$ by $\sigma$ as well. Namely, $\sigma$ acting on $R_d$ is the Frobenius automorphism of $R_d$ which fixes $\pi$ such that the reduction map $\sigma$ modulo $\pi$ is the $q$-th power map of $\mathbf{F}_{q^d}$.

Via the map $x: A \to R_d$, one obtains the fibre $M_x$ at $x$:

$$M_x = M \otimes_A R_d$$

which is a free module of finite rank over $R_d$. The fibre map (evaluating $\phi$ at $x$)

$$\phi_x = \phi \otimes_A R_d : M_x \longrightarrow M_x$$

is $\sigma$-linear (not linear in general) over $R_d$. But its $d$-th iterate $\phi_x^d$ is $R_d$-linear. Thus, the characteristic polynomial $\det(I - T\phi_x^d)$ of $\phi_x^d$ acting on the free $R_d$-module $M_x$ is well defined. By working with a basis, one can easily show that this characteristic polynomial actually has coefficients in $R$, not just in $R_d$. Alternatively, the map $\phi_x^d$ can be defined over $R$. The characteristic polynomial is independent of the choice of the geometric point $\bar{x}$ of the closed point on $X/\mathbf{F}_q$ containing $\bar{x}$. Thus, we can speak of the characteristic polynomial at a closed point $\bar{x}$ on $X/\mathbf{F}_q$.

**Definition 2.3.** The L-function of a $\sigma$-module $(M, \phi)$ over $X/\mathbf{F}_q$ is defined by the Euler product

$$(2.1) \qquad L(\phi, T) = \prod_{\bar{x} \in X/\mathbf{F}_q} \frac{1}{\det(I - \phi_x^{\deg(x)} T^{\deg(x)})} \in R[[T]],$$

where $\bar{x}$ runs over the closed points of $X/\mathbf{F}_q$ and $x$ is the Teichmüller lifting of $\bar{x}$ defined with respect to $\sigma$. For simplicity, we have denoted $\deg(\bar{x})$ by $\deg(x)$.

For example, if $\phi$ is the identity $\sigma$-module $(A_0, \sigma)$, then $L(\phi, T)$ is just the zeta function of the affine variety $X/\mathbf{F}_q$ which is a rational function by Dwork's rationality theorem. The above definition shows that the L-function $L(\phi, T)$ is always a power series with coefficients in $R$. Thus, it is trivially analytic in the open unit disk $|T|_\pi < 1$. We shall normalize the $p$-adic absolute value by requiring $|\pi|_\pi = 1/p$. It is true but non-trivial that the L-function is meromorphic on the closed unit disk $|T|_\pi \leq 1$. It is shown in [18] that the L-function $L(\phi, T)$ in (2.1) is in general not a meromorphic function. Thus, to get further meromorphic continuation of $L(\phi, T)$, we must impose additional conditions on $(M, \phi)$. The generalization of the Monsky trace formula as described in section 3 shows that $L(\phi, T)$ is $p$-adic meromorphic if $\phi$ is overconvergent.

We now describe Dwork's conjecture in the more general setting of $\sigma$-modules. Let $(M, \phi)$ be a $\sigma$-module over $X/\mathbf{F}_q$. For a rational number $s \in \mathbf{Q}$, we define the slope $s$ part of the characteristic polynomial of $(M, \phi)$ at $x$ to be the partial product

$$\det{}_s(\phi_x, T) = \prod_{\mathrm{ord}_{\pi^{\deg(x)}}(\alpha) = s} (1 - \alpha T),$$

where $\alpha$ runs over the reciprocal roots of slope $s$ of the characteristic polynomial $\det(I - \phi_x^{\deg(x)} T)$ defined with respect to the valuation $\mathrm{ord}_{\pi^{\deg(x)}}$. For each $s \in \mathbf{Q}$, the polynomial $\det_s(\phi_x, T)$ actually has coefficients in $R$ although its roots are in general in an extension field. In our setting, the factor $\det_s(\phi_x, T)$ is non-trivial only when $s$ is a non-negative rational number and thus $\det_s(\phi_x, T) = 1$ for negative number $s$. It is clear that we have the purity or slope decomposition

(2.2) $$\det(I - \phi_x^{\deg(x)} T) = \prod_{s \in \mathbf{Q}} \det{}_s(\phi_x, T).$$

For each $x$, this is a finite product since a polynomial has only finitely many zeros. The slope $s$ zeta function of $(M, \phi)$ is defined to be the following Euler product

$$L_s(\phi, T) = \prod_{\bar{x} \in X/\mathbf{F}_q} \frac{1}{\det{}_s(\phi_x, T^{\deg(x)})} \in R[[\pi^s T]],$$

where $\bar{x}$ runs over the closed points of $X/\mathbf{F}_q$. Strictly speaking, we need to assume that $\pi^s \in R$, but this can always be achieved by going to a totally ramified finite extension of $R$ if necessary. Clearly, $L_s(\phi, T)$ is analytic in $|T|_\pi < p^s$. By (2.2), we have the decomposition formula for $L(\phi, T)$:

(2.3) $$L(\phi, T) = \prod_{s \in \mathbf{Q}} L_s(\phi, T).$$

For an arbitrary $\sigma$-module $(M, \phi)$, one cannot expect that the slope $s$ zeta function $L_s(\phi, T)$ to be meromorphic since their product over all $s$ may already fail to be meromorphic. However, in geometric situation, the $\sigma$-module $(M, \phi)$ is overconvergent and thus the total L-function $L(\phi, T)$ is known to be meromorphic. This leads to Dwork's conjecture.

**Conjecture 2.4** (Dwork). *Let $(M, \phi)$ be an overconvergent $\sigma$-module over $X/\mathbf{F}_q$. Then for each rational number $s \in \mathbf{Q}$, the slope $s$ zeta function $L_s(\phi, T)$ is $p$-adic meromorphic.*

In order to prove the conjecture, we have to understand how the slope decomposition in (2.2) varies as $x$ varies. Namely, we need to understand how the Euler factor of $L_s(\phi, T)$ at $x$ varies as $x$ varies. The first step is to understand how the decomposition pattern (Newton polygon) in (2.2) varies as $x$ varies. Grothendieck's specialization theorem says that the decomposition pattern in (2.2) varies in an algebraic way with $x$. In particular, this implies that there are only finitely many possible decomposition patterns for a given $(M, \phi)$ and thus the product in (2.3) is actually a finite product. For a fixed decomposition pattern and a fixed slope $s$, we would like to pull out the slope $s$ part $\phi_s$ of $\phi$ in an algebraic way so that $\phi_s$ also forms a $\sigma$-module over $X$. This is indeed the case by the Hodge-Newton decomposition if $(M, \phi)$ is ordinary fibre by fibre, namely, if the Newton polygon coincides with the Hodge polygon fibre by fibre. In the harder non-ordinary case, Katz's isogeny theorem shows that by shrinking $X$ if necessary, we can still pull out the slope $s$ part $\phi_s$ up to isogeny such that $\phi_s$ is a $\sigma$-module over $X$. This at least reduces the slope $s$ zeta function $L_s(\phi, T)$ to the L-function $L(\phi_s, T)$ of the more natural $\sigma$-module $\phi_s$ where one could hope for some tools such as a trace formula. If it happens that the $\sigma$-module $\phi_s$ is also overconvergent, then we are done by Monsky's trace formula. Since the Hodge-Newton decomposition does not exist over $A$ (it exists only over the completion $A_0$ of $A$), these pure slope $\sigma$-modules $\phi_s$ are NOT overconvergent in general. They rarely are in fact. Thus, the existing overconvergent theory does not suffice and we have to introduce more new ideas.

The aim of the present paper is to introduce a reduction approach to reduce higher rank $\phi_s$ to rank one $\phi_s$ so that our previous limiting approach can be applied. The decomposition formula in (2.3) was the starting point of our limiting approach. It is also the starting point for the reduction approach. The full generality of our result says that if $\psi$ is in the tensor category generated by the $\phi_s$ as $s$ and the overconvergent $\phi$ vary, then the L-function $L(\psi, T)$ is meromorphic.

## 3. The Monsky trace formula

As mentioned above, if a $\sigma$-module $\phi$ is overconvergent, the Monsky trace formula can be generalized to prove the meromorphic continuation of the L-function $L(\phi, T)$. This trace formula will be used as a tool in our work. We briefly recall the general trace formula and the related construction in the form we need. For a proof, see the appendix in section 10. Throughout this section, we assume that $(M, \phi)$ is overconvergent. Recall that our ring $A$ is an $n$-dimensional overconvergent integral domain. We do not shrink $X$ yet in this section. Thus, we allow projective $\sigma$-modules in this section.

First, we restrict to the identity overconvergent $\sigma$-module with $(M, \phi) = (A, \sigma)$. For a non-negative integer $i$, let $\Omega^i A$ be the $A$-module of differential $i$-forms on $A$ over $R$. Since $X$ is smooth and $A$ is flat, $\Omega^i A$ is a finite projective $A$-module. In particular, $\Omega^n A$ is projective of rank 1 over $A$. The Frobenius map $\sigma$ extends to $\Omega^i A$ as an injective $R$-linear endomorphism, denoted by

$$\sigma_i : \Omega^i A \longrightarrow \Omega^i A.$$

Clearly,
$$\sigma_1(b) \equiv 0 \ (\mathrm{mod}\, \pi), \ b \in \Omega^1 A.$$

Since $\sigma_i = \wedge^i \sigma_1$, we deduce

(3.1) $$\sigma_i(b) \equiv 0 \ (\mathrm{mod}\, \pi^i), \ b \in \Omega^i A.$$

There is a trace map
$$\mathrm{Tr}_i : \Omega^i A \longrightarrow \sigma_i(\Omega^i A)$$

such that for $a \in A, b \in \Omega^i A$, we have
$$\mathrm{Tr}_i(\sigma(a)b) = \sigma(a)\mathrm{Tr}_i(b).$$

Furthermore, as endomorphisms on $\Omega^i A$,

(3.2) $$\sigma_i^{-1} \circ \mathrm{Tr}_i \circ \sigma_i = [\bar{A} : \sigma(\bar{A})] = q^n,$$

where $\bar{A}$ means the reduction modulo $\pi$ of $A$. Thus, the map $\sigma_i$ acting on $\Omega^i A$ has a one-sided left inverse when tensored with $\mathbf{Q}$. This one sided inverse is the map $\theta_i$ defined by
$$\theta_i = \sigma_i^{-1} \circ \mathrm{Tr}_i : \ \Omega^i A \longrightarrow \Omega^i A.$$

The map $\theta_i$ is an example of a Dwork operator, which by definition is $\sigma^{-1}$-linear in the following sense:
$$\theta_i(\sigma(a)m) = a\theta_i(m), \ a \in A, m \in \Omega^i A.$$

Dwork operators are of trace class and play an important role.

For an overconvergent $\sigma$-module $(M, \phi)$, its underlying dual $A$-module $M^\vee$ is defined by
$$M^\vee = \mathrm{Hom}_A(M, A).$$

Define a dual map $\phi^\vee : M^\vee \to M^{\sigma\vee}$ by
$$\phi^\vee(f)(m) = f(\phi(m)), \ f \in M^\vee, m \in M.$$

The pair $(M^\vee, \phi^\vee)$ is not a $\sigma$-module yet, because the arrow direction of the map $\phi^\vee$ is not the right one. In the case that $\phi$ is a unit root $\sigma$-module, the contragradient (the inverse of $\phi^\vee$) does give an overconvergent $\sigma$-module structure on $M^\vee$. We shall not need this construction. To describe the trace formula, we need to discuss a different type of "duality" which turns a $\sigma$-module to a module with a Dwork operator.

Let $(M, \phi)$ be an overconvergent $\sigma$-module. We define a new $A$-module by
$$M^* = \mathrm{Hom}_A(M, \Omega^n A).$$

Define a Dwork operator $\phi^*$ on $M^*$ as follow: If $f \in M^*$ and $m \in M$, then
$$\phi^*(f)(m) = (\sigma_n^{-1} \circ \mathrm{Tr}_n)(f(\phi(m))) = \theta_n(f(\phi(m))).$$

One checks that $\phi^*(f)$ is in $M^*$. Furthermore, the map $\phi^*$ is $\sigma^{-1}$-linear. That is, for $a \in A$, we have
$$\phi^*(\sigma(a)f) = a\phi^*(f).$$

For each non-negative integer $0 \leq i \leq n$, let
$$\Omega^i M = M \otimes_A \Omega^i A.$$

This is again a finite projective $A$-module. Let
$$\phi_i = \phi \otimes \sigma_i,$$

where $(\Omega^i A, \sigma_i)$ is the $\sigma$-module given by the action of $\sigma_i$ on $\Omega^i A$. Then, the pair $(\Omega^i M, \phi_i)$ becomes an overconvergent $\sigma$-module (the projective version of a $\sigma$-module). By (3.1), we obtain the congruence

(3.3) $$\phi_i \equiv 0 \ (\mathrm{mod} \pi^i).$$

Let
$$M_i^* = \mathrm{Hom}_A(\Omega^i M, \Omega^n A).$$

In a similar way, we define the corresponding Dwork operator $\phi_i^*$ on $M_i^*$ by
$$\phi_i^*(f)(m) = (\sigma_n^{-1} \circ \mathrm{Tr}_n)(f(\phi_i(m))) = \theta_n(f(\phi_i(m))), \ m \in \Omega^i M, \ f \in M_i^*.$$
This is a $\sigma^{-1}$-linear operator. For $i = 0$, we have
$$(\Omega^0 A, \sigma_0) = (A, \sigma), \ (M_0^*, \phi_0^*) = (M^*, \phi^*).$$
The congruence in (3.3) shows for all $i$,

(3.4) $$\phi_i^* \equiv 0 \ (\mathrm{mod} \pi^i).$$

Let $K$ be the quotient field of $R$. The finite $A \otimes K$-module $M_i^* \otimes_R K$ is an infinite dimensional vector space over $K$, which is a direct limit of a sequence of $p$-adic Banach spaces with orthonormal bases. Since $\phi$ is overconvergent, the $K$-linear Dwork operator $\phi_i^*$ is a nuclear operator of the $p$-adic space $M_i^* \otimes_R K$ over $K$. It follows that the Fredholm determinant $\det(I - T\phi_i^*|M_i^* \otimes K)$ is well defined and is a $p$-adic entire function with coefficients in $K$. The coefficients are actually in $R$ since the above construction shows that $\phi_i$ is integral. The generalized Monsky trace formula can then be stated as

**Theorem 3.1.** *Let $(M, \phi)$ be an overconvergent $\sigma$-module over $X/\mathbf{F}_q$. Then*

(3.5) $$L(\phi, T) = \prod_{i=0}^{n} \det(I - T\phi_{n-i}^*|M_{n-i}^* \otimes_R K)^{(-1)^{i-1}}.$$

If $M$ has rank one, Theorem 3.1 is an immediate consequence of Theorem 5.3 in [13]. The general case is given in the appendix of this paper, see Theorem 10.10. In particular, we obtain

**Corollary 3.2.** *Let $(M, \phi)$ be an overconvergent $\sigma$-module over $X/\mathbf{F}_q$. Then the L-function $L(\phi, T)$ is p-adic meromorphic.*

Next, we give an alternative description of the formula in (3.5) which is perhaps more familiar. The perfect pairing
$$\Omega^i A \times \Omega^{n-i} A \longrightarrow \Omega^n A : \ (w_1, w_2) \longrightarrow w_1 \wedge w_2$$
gives a natural identification
$$\Omega^i A = \mathrm{Hom}_A(\Omega^{n-i}A, \Omega^n A) = A_{n-i}^*.$$
Thus, for the identity $\sigma$-module $(A, \sigma)$, the $(n-i)$-th Dwork operator $\sigma_{n-i}^*$ acting on $A_{n-i}^*$ becomes the Dwork operator $\theta_i$ acting on $\Omega^i A$:
$$\sigma_{n-i}^* = \theta_i.$$
In fact, one checks that for $w_1 \in \Omega^i A$ and $w_2 \in \Omega^{n-i} A$,
$$\begin{aligned}
\sigma_{n-i}^*(w_1)(w_2) &= \sigma_n^{-1} \circ \mathrm{Tr}_n(w_1 \wedge \sigma_{n-i}(w_2)) \\
&= \sigma_n^{-1}(\mathrm{Tr}_i(w_1) \wedge \sigma_{n-i}(w_2)) \\
&= (\sigma_i^{-1} \circ \mathrm{Tr}_i)(w_1) \wedge w_2 \\
&= \theta_i(w_1) \wedge w_2.
\end{aligned}$$
More generally, for an overconvergent $\sigma$-module $(M, \phi)$, we can write
$$\begin{aligned}
M_i^* &= \mathrm{Hom}_A(M \otimes_A \Omega^i A, \Omega^n A) \\
&= \mathrm{Hom}_A(M, A) \otimes_A \mathrm{Hom}_A(\Omega^i A, \Omega^n A) \\
&= M^\vee \otimes_A \Omega^{n-i} A \\
&= \Omega^{n-i} M^\vee.
\end{aligned}$$

Under this identification, the Dwork operator $\phi_i^*$ acting on $M_i^*$ becomes a Dwork operator on $\Omega^{n-i} M^\vee$, denoted by $\theta_{n-i}(\phi)$:

$$\theta_{n-i}(\phi) : \Omega^{n-i} M^\vee \longrightarrow \Omega^{n-i} M^\vee.$$

In these new notations, the trace formula in (3.5) can then be written as

$$(3.6) \qquad L(\phi, T) = \prod_{i=0}^{n} \det(I - T\theta_i(\phi) | \Omega^i M^\vee \otimes_R K)^{(-1)^{i-1}}.$$

If $(M, \phi)$ is the identity overconvergent $\sigma$-module $(A, \sigma)$, one can refine the chain level formula in (3.6) to get a cohomological formula. In this case, the $\sigma$-module $(A, \sigma)$ is self-dual. Thus, the Dwork operator $\theta_i(\sigma)$ on $\Omega^i A^\vee = \Omega^i A$ is given by

$$\theta_i(\sigma) = \theta_i = \sigma_i^{-1} \circ \text{Tr}_i.$$

These operators $\theta_i$ commute with the exterior differentiation $d$ of the De Rham complex of $A \otimes K$:

$$\Omega^\cdot A \otimes K : \quad 0 \longrightarrow A \otimes K \xrightarrow{d} \Omega^1 A \otimes K \xrightarrow{d} \Omega^2 A \otimes K \xrightarrow{d} \cdots.$$

Namely,

$$\theta_i \circ d = d \circ \theta_{i-1}.$$

The cohomology of the De Rham complex $\Omega^\cdot A \otimes K$ is the formal cohomology $H^i_{\text{MW}}(X, K)$ of Monsky-Washnitzer [13]. This is now known to be a finite dimensional vector space over $K$, by the work of Mebkhout [12] and Berthelot [2]. Since the $\theta_i$ commutes with the exterior differentiation $d$, the Dwork operator induces a chain map of the De Rham complex. Serre's $p$-adic spectral theory [16] implies that the chain level formula in (3.6) passes to the cohomology level. It follows that the zeta function $Z(X, T)$ of $X$ over $\mathbf{F}_q$ (the L-function of the identity $\sigma$-module) can be computed using the Monsky-Washnitzer cohomology:

$$Z(X, T) = \prod_{i=0}^{n} \det(I - T\theta_i | H^i_{\text{MW}}(X, K))^{(-1)^{i-1}}.$$

Since $\theta_i \circ \sigma_i = q^n$, the above formula can be rewritten as

$$Z(X, T) = \prod_{i=0}^{n} \det(I - Tq^n \sigma_i^{-1} | H^i_{\text{MW}}(X, K))^{(-1)^{i-1}}.$$

More generally, assume that our given overconvergent $\sigma$-module $(M, \phi)$ is endowed with an integrable overconvergent connection $\nabla$ such that $(M \otimes K, \phi, \nabla)$ becomes an overconvergent $F$-isocrystal [1] and $\phi$ is the Frobenius map $\sigma$ acting on $M$. The connection $\nabla$ gives rise to the De Rham complex of $M \otimes K$:

$$\Omega^\cdot M \otimes K : \quad 0 \longrightarrow M \otimes K \xrightarrow{\nabla} \Omega^1 M \otimes K \xrightarrow{\nabla} \Omega^2 M \otimes K \xrightarrow{\nabla} \cdots.$$

The cohomology $H^i(\Omega^\cdot M \otimes K)$ of this complex is the rigid cohomology $H^i_{\text{rig}}(X, M \otimes K)$ of Berthelot [1], which is conjectured to be a finite dimensional vector space over $K$. As shown in Etesse-LeStum [7], the L-function $L(\phi, T)$ in this case can be similarly expressed in terms of the rigid cohomology $H^i_{\text{rig}}(X, M^\vee)$ of the dual overconvergent $F$-isocrystal $M^\vee$:

$$L(\phi, T) = \prod_{i=0}^{n} \det(I - Tq^n \sigma_i^{-1} | H^i_{\text{rig}}(X, M^\vee))^{(-1)^{i-1}},$$

where $q^n \sigma_i^{-1}$ acts as a nuclear operator on $H^i_{\text{rig}}(X, M^\vee)$. If one uses the rigid cohomology $H^\cdot_{\text{rig},c}$ with compact support, the cohomological formula becomes

$$L(\phi, T) = \prod_{i=0}^{2n} \det(I - T\sigma_i | H^i_{\text{rig},c}(X, M))^{(-1)^{i-1}},$$

where $H^i_{\text{rig},c}(X, M) = 0$ for $0 \leq i < n$ since our space $X$ is affine. The L-function $L(\phi, T)$ in such an overconvergent $F$-isocrystal situation is conjectured to be a rational function. This follows from the conjectural finite dimensionality of the rigid cohomology, which in turns follows from the semi-stable reduction conjecture for overconvergent F-crystals, see A.J. de Jong [4] for more detail.

In this paper, we study the more general category of $\sigma$-modules with no overconvergent connections and hence no cohomological formulas involved. The L-functions are not rational and in fact not even meromorphic in such a generality.

To conclude this section, we point out a weak entireness result about the L-function. This result is not needed for our proof of Dwork's conjecture. But it is an interesting consequence of the above description. In fact, the congruence in (3.4) shows that

$$\det(I - \phi_i^* T) \in 1 + \pi^i T R[[\pi^i T]].$$

Combining this with Theorem 3.1, we obtain

**Corollary 3.3.** *Let $(M, \phi)$ be an overconvergent $\sigma$-module over $X/\mathbf{F}_q$. Then*

$$(3.7) \qquad \frac{L(\phi, T)^{(-1)^{n-1}}}{\det(I - T\phi^* | M^* \otimes_R K)} \in 1 + \pi T R[[\pi T]].$$

This result shows that the quotient in (3.7) has no zeros and poles in the disc $\text{ord}_\pi(T) > -1$ (the disc $|T|_\pi < p$). Alternatively, the L-function $L(\phi, T)^{(-1)^{n-1}}$ has no poles in the disc $|T|_\pi < p$.

**Corollary 3.4.** *Let $X$ be a connected smooth affine variety over $\mathbf{F}_q$ of dimension $n$. Then, the zeta function of $X/\mathbf{F}_q$ has the form*

$$(3.8) \qquad Z(X, T)^{(-1)^{n-1}} = \frac{P(T)}{\prod_{j=1}^m (1 - q\alpha_j T)},$$

*where $P(T)$ is a polynomial with integer coefficients and the $\alpha_j$ are p-adic integers.*

The proof is similar. It suffices to take $(M, \phi) = (A, \sigma)$ in the above. The reason that we can use $q$ instead of $\pi$ in (3.8) is that there is a Frobenius lifting $\sigma$ on some overconvergent ring $A$ lifting $X$ such that

$$\sigma_i \equiv 0 \pmod{q^i}, \ \phi_i \equiv 0 \pmod{q^i}, \ \phi_i^* \equiv 0 \pmod{q^i}.$$

A little cohomological argument shows that Corollary 3.4 remains true if $X$ is more generally a connected $n$-dimensional set theoretic complete intersection contained in some connected smooth affine variety over $\mathbf{F}_q$. Such a result can also be proved using $\ell$-adic cohomology and Deligne's integrality theorem, see [22] for some further information on such weak entireness results. In particular, the shifted reciprocal poles $\alpha_j$ in Corollary 3.4 are algebraic integers not just $p$-adic integers.

## 4. Hodge-Newton decomposition

In this section, we define the basis sequence, basis filtration and basis polygon of a $\sigma$-module under a given basis. For our purpose, the basis polygon is simpler to define and more convenient to use than the Hodge polygon. Our basis polygon is a lower bound for both the Newton polygon and the Hodge polygon fibre by fibre. We shall review Grothendieck's specialization theorem for Newton polygons and the Hodge-Newton decomposition as proved by Dwork [5][6] and more generally by Katz [10]. This allows us to reformulate Dwork's conjecture in the ordinary case in terms of L-functions of pure slope $\sigma$-modules.

**Definition 4.1.** Let $(M, \phi)$ be a $\sigma$-module with a given basis $\vec{e} = \{e_1, e_2, \cdots\}$ over $A_0$. For each integer $i \geq 1$, let $d_i$ be the smallest positive integer $d$ such that for all $j > d$, we have
$$\phi(e_j) \equiv 0 \pmod{\pi^i}.$$
For $i \geq 0$, we define
$$h_i = d_{i+1} - d_i, \ d_0 = 0.$$
The sequence $h = h(\vec{e}) = \{h_0, h_1, \cdots\}$ is called the **basis sequence** of $\phi$ defined with respect to the basis $\vec{e}$. If we let $M_{(i)}$ be the finite free $A_0$-module with the basis $\{e_{d_i+1}, e_{d_i+2}, \cdots\}$, then we have a decreasing **basis filtration** of finite free $A_0$-submodules:
$$M = M_{(0)} \supset M_{(1)} \supset \cdots \supset M_{(j)} \supset \cdots \supset 0, \ \cap_j M_{(j)} = 0,$$
where each quotient $M_{(i)}/M_{(i+1)}$ is also a finite free $A_0$-module and
$$\phi(M_{(i)}) \subseteq \pi^i M, \ \text{rank}(M_{(i)}/M_{(i+1)}) = h_i.$$
Equivalently, the matrix $G$ of $\phi$ with respect to $\vec{e}$ (defined by $\phi(\vec{e}) = \vec{e}G$) is of the form
$$G = (G_0, \pi G_1, \pi^2 G_2, \cdots, \pi^i G_i, \cdots),$$
where each block $G_i$ is a matrix over $A_0$ with $h_i$ columns.

To describe Dwork's conjecture in the ordinary case, we need to define the notion of basis polygon and Newton polygon.

**Definition 4.2.** Let $(M, \phi)$ be a $\sigma$-module with a basis $\vec{e}$. Let $h = h(\vec{e}) = \{h_0, h_1, \cdots, h_k\}$ be the basis sequence defined with respect to $\vec{e}$. We define the **basis polygon** $P(\vec{e})$ of $(M, \phi)$ with respect to $\vec{e}$ to be the convex closure in the plane of the following lattice points:
$$(0,0), (h_0, 0), (h_0 + h_1, h_1), (h_0 + h_1 + h_2, h_1 + 2h_2), \cdots$$
$$(h_0 + \cdots + h_k, h_1 + 2h_2 + \cdots + kh_k).$$
Namely, the basis polygon is the polygon with a side of slope $i$ and horizontal length $h_i$ for every integer $0 \leq i \leq k$.

It is clear that our definition of the basis polygon $P(\vec{e})$ and the basis sequence $h(\vec{e})$ depends on the given basis $\vec{e}$. The basis polygon and the basis sequence are somewhat similar to Hodge polygon and Hodge numbers. But they are different in general. The same $\sigma$-module with different bases will give rise to different basis polygons and different basis sequences. The abstract Hodge polygon as discussed in [10] can only be defined generically or at a point. One can show that the generic Hodge polygon (or the Hodge polygon at each fibre) lies on or above our basis polygon. In the ordinary case to be defined below, the basis $\vec{e}$ will be taken to be the optimal one and our basis polygon is the same as the Hodge polygon. Since we are only concerned with the ordinary case in this section, nothing is lost. In the more general case, we work directly with the Newton polygon. The Hodge polygon and basis polygon play only an auxiliary role. It is the Newton polygon that is most important for us.

**Definition 4.3.** Let $(M, \phi)$ be a $\sigma$-module over $X/\mathbf{F}_q$. For each closed point $\bar{x} \in X/\mathbf{F}_q$, the Newton polygon of a $\sigma$-module $(M, \phi)$ at $\bar{x}$ is the Newton polygon of the characteristic polynomial $\det(I - \phi_x^{\deg(x)} T)$ defined with respect to the valuation $\mathrm{ord}_{\pi^{\deg(x)}}$.

The Newton polygon has several other interpretations, see [10] for more details and for various notions related to Newton polygons. We note that the notation $A_0$ in [10] corresponds to the reduction $A_0/\pi A_0$ in this paper.

The abstract version of Mazur's theorem [11] says that the Newton polygon lies above the Hodge polygon fibre by fibre. Since the basis polygon $P(\vec{e})$ is a lower bound for the Hodge polygon fibre by fibre, we deduce the following result.

**Lemma 4.4.** *Let $(M, \phi)$ be a $\sigma$-module over $X/\mathbf{F}_q$ with a basis $\vec{e}$. The Newton polygon of $(M, \phi)$ at each fibre lies on or above the basis polygon $P(\vec{e})$ of $(M, \phi)$.*

This can also be proved directly since our definition of the basis filtration immediately implies that

$$\mathrm{ord}_\pi(\wedge^{h_0 + \cdots + h_i} \phi) \geq h_1 + 2h_2 + \cdots + ih_i.$$

As for the Newton polygon, Grothendieck's specialization theorem holds for $\sigma$-modules as well. Instead of reproving this type of results, we shall take the shorter approach of reducing them to the familiar case of F-crystals as treated in [10]. One difference between a $\sigma$-module and an F-crystal is that the map $\phi$ of a $\sigma$-module may not be injective. This is the only difference if the base space $X$ is a point. For many purpose, the situation can be reduced to the injective case. The following lemma holds for both convergent and overconvergent $\sigma$-modules.

**Lemma 4.5.** *Let $(M, \phi)$ be a $\sigma$-module (resp. an overconvergent $\sigma$-module) over $X$. Then, shrinking $X$ if necessary, there is a $\phi$-stable filtration of free $A_0$-modules (resp. free $A$-modules) with free quotients:*

$$0 \subset N_0 \subset M$$

*such that the restriction of $\phi$ to $N_0$ is nilpotent and the quotient $(M/N_0, \phi)$ is an injective $\sigma$-module.*

*Proof.* Let us work with $A_0$-modules since the proof for $A$-modules is completely similar. Let $L$ be the localization of $A_0$ at the prime ideal $\pi$. That is, we invert all elements in $A_0$ which are not divisible by $\pi$. Let

$$N_0 = \bigcup_{k=1}^{\infty} \mathrm{Ker}(\phi^k | M).$$

This is the nilpotent part of $M$ for the $\phi$-action. It is clear that both $N_0$ and $M/N_0$ are torsion free. They are both finitely generated since $A_0$ is Noetherian. The $L$-module $N_0 \otimes_{A_0} L$ is the nilpotent part of $M \otimes_{A_0} L$ for the $\phi$-action. Since $L$ is a principal ideal domain, $N_0 \otimes_{A_0} L$ is a free $L$-submodule of $M \otimes_{A_0} L$ whose quotient is torsion free and hence free. Thus, we have a $\phi$-stable filtration of finite free $L$-modules with free quotients

$$0 \subset N_0 \otimes_{A_0} L \subset M \otimes_{A_0} L.$$

Let $\{e_1, \cdots, e_i, e_{i+1}, \cdots, e_r\}$ be a basis of the free $L$-module $M \otimes_{A_0} L$ such that $\{e_1, \cdots, e_i\}$ is a basis of the submodule $N_0 \otimes_{A_0} L$. Shrinking $X$ if necessary, we see that $\{e_1, \cdots, e_r\}$ is a basis of the free $A_0$-module $M$ such that $\{e_1, \cdots, e_i\}$ is a basis of the submodule $N_0$. The proof is complete. □

The nilpotent part $N_0$ has no contributions at all for the L-function. Thus, Lemma 4.5 shows that for L-function purpose, we can always restrict our attention to injective $\sigma$-modules. Shrinking a little further, we may even assume that $\phi$ is injective fibre by fibre. At a point, our notion of an injective $\sigma$-module coincides with the notion of an F-crystal. Thus, all pointwise results on F-crystals in [10] carry over to injective $\sigma$-modules at a point. These pointwise results are the key to all the global results in [10]. Repeating the proofs in [10], we conclude that all global results in [10] carry over to $\sigma$-modules which are fibre by fibre injective, except that we may need to shrink $X$ a little further in order to replace a locally free module (which is the setup of [10]) by a free module (which is our setup).

In particular, we have

**Lemma 4.6.** *Let $(M, \phi)$ be a $\sigma$-module of rank $r$, which is fibre by fibre injective. Let $s$ be a non-negative rational number. Then, all Newton slopes of $\phi$ at all fibres are at least $s$ if and only if for every positive integer $k$, we have the inequality*

$$\mathrm{ord}_\pi \phi^{k+r-1} \geq ks.$$

This result holds fibre by fibre by the basic slope estimate in [10] and thus it holds globally.

**Theorem 4.7** (Grothendieck). *A $\sigma$-module $(M, \phi)$ over $X/\mathbf{F}_q$ has a generic Newton polygon. Furthermore, the Newton polygon of $(M, \phi)$ at each fibre lies on or above the generic Newton polygon.*

*Proof.* We give an outline of the proof following [10]. By Lemma 4.5, we may assume that $\phi$ is fibre by fibre injective. By the Dieudonne-Manin theorem for an F-crystal at a point, each Newton slope of $\phi$ at each fibre is a non-negative rational number with denominator dividing $r!$, where $r$ is the rank of $\phi$. Applying various exterior power constructions of $\phi$, it suffices to prove the claim that for each non-negative rational number $s$, the set of points in $X$ at which all Newton slopes of $\phi$ are at least $s$ is Zariski closed. Fix a large positive integer $m$ depending on $s$ such that

$$(4.1) \qquad s - \frac{1}{r!} < \frac{m}{m+r-1} s < s.$$

Lemma 4.6 shows that if $\phi$ has all Newton slopes at least $s$, then

$$(4.2) \qquad \mathrm{ord}_\pi \phi^{m+r-1} \geq ms.$$

Conversely, if inequality (4.2) holds, then clearly all Newton slopes of $\phi$ are at least $sm/(m+r-1)$. This together with (4.1) imply that $\phi$ has all Newton slopes at least $s$. Thus, $\phi$ has all Newton slopes at least $s$ if and only if (4.2) holds. But the set of points in $X$ for which (4.2) holds is clearly a Zariski closed subset of $X$. The proof is complete. □

**Definition 4.8.** Let $\vec{e}$ be a basis of a $\sigma$-module $(M, \phi)$. The basis $\vec{e}$ is called **ordinary** if the basis polygon $P(\vec{e})$ of $(M, \phi)$ coincides with the Newton polygon of each fibre $(M, \phi)_x$, where $x$ is the Teichmüller lifting of the closed point $\bar{x} \in X/\mathbf{F}_q$. The basis filtration attached to $\vec{e}$ is called ordinary if the basis $\vec{e}$ is ordinary. Similarly, the basis $\vec{e}$ is called ordinary up to slope $j$ side (resp. ordinary up to slope $j$) if the basis polygon $P(\vec{e})$ coincides with the Newton polygon of each fibre for all sides up to slope $j$ side including the whole slope $j$ side (resp. all sides up to slope $j$ including only a non-trivial portion of the slope $j$ side). In particular, the basis $\vec{e}$ is said to be ordinary at the slope zero side (resp. at slope zero) if the Newton polygon at each fibre has a horizontal side of length $h_0(\vec{e})$ (resp. a horizontal side of positive length if $h_0 > 0$). We say that $(M, \phi)$ is ordinary (resp. ordinary up to slope $j$ side, resp. ordinary up to slope $j$) if it has a basis which is ordinary (resp. ordinary up to slope $j$ side, resp. ordinary up to slope $j$).

One checks that the category of ordinary $\sigma$-modules over $X$ is closed under direct sum, tensor product, symmetric power and exterior power. For example, if $\phi$ is a $\sigma$-module with a basis sequence $h = (h_0, h_1, \cdots, h_k)$ and $\psi$ is a $\sigma$-module with a basis sequence $g = (g_0, g_1, \cdots, g_\ell)$. Then the tensor product $\phi \otimes \psi$ is a $\sigma$-module with a basis sequence

$$h \otimes g = (h_0 g_0, h_0 g_1 + h_1 g_0, h_0 g_2 + h_1 g_1 + h_2 g_0, \cdots, h_k g_\ell).$$

If both $\phi$ and $\psi$ are ordinary, one checks that $\phi \otimes \psi$ is ordinary. Similarly, the exterior product $(\wedge^2 M, \wedge^2 \phi)$ is a $\sigma$-module with a basis sequence

$$\wedge^2 h = (\frac{h_0^2 - h_0}{2}, h_0 h_1, h_0 h_2 + \frac{h_1^2 - h_1}{2}, h_0 h_3 + h_1 h_2, h_0 h_4 + h_1 h_3 + \frac{h_2^2 - h_2}{2}, \cdots).$$

If $\phi$ is ordinary, then $\wedge^2 \phi$ is ordinary. The following ordinary Hodge-Newton decomposition is due to Dwork (Lemma 5.1 in [6]), see also Lemma 3.4 in [19] for a quick proof.

**Theorem 4.9.** Let $(M, \phi)$ be a $\sigma$-module ordinary up to slope $j$ side, where $0 \leq j \leq k$ is a non-negative integer. Then, there is an increasing $\phi$-stable filtration of free $A_0$-submodules of finite rank,

$$0 \subset M_0 \subset M_1 \subset \cdots \subset M_j \subset M$$

which is transversal to the ordinary basis filtration:

$$M = M_i \oplus M_{(i+1)}, \ 0 \leq i \leq j.$$

Furthermore, for $0 \leq i \leq j$, the quotient is of the form

$$(M_i/M_{i-1}, \phi) = (U_i, \pi^i \phi_i),$$

where $(U_i, \phi_i)$ is a unit root $\sigma$-module of rank $h_i$.

We shall call $(U_i, \phi_i)$ the unit root $\sigma$-module coming from the slope $i$ part of $(M, \phi)$. It is not hard to see that if $(M, \phi)$ is ordinary up to slope $i$ side, then we have

(4.3) $\qquad \text{ord}_\pi(\wedge^{h_0 + \cdots + h_i} \phi) = h_1 + 2h_2 + \cdots + ih_i.$

Note that the Hodge-Newton decomposition is false if we replace $A_0$ by $A$. Namely, if we start with an **overconvergent** and ordinary $\sigma$-module $\phi$, the Hodge-Newton decomposition cannot be obtained over $A$. It exists only over the completion $A_0$ of $A$. The unit root pieces $\phi_i$ obtained from the Hodge-Newton decomposition cannot be defined over $A$ and thus will NOT be overconvergent in general, even if the original ambient $\phi$ is overconvergent. Using the notation of (2.3), we have

$$L_i(\phi, T) = L(\phi_i, \pi^i T).$$

Thus, Conjecture 2.3 in the ordinary case is equivalent to the following

**Conjecture 4.10** (Dwork). *Let $(M, \phi)$ be an **overconvergent** ordinary $\sigma$-module over $X/\mathbf{F}_q$. Let $(U_i, \phi_i)$ be the unit root $\sigma$-module coming from the slope $i$ part of $(M, \phi)$, where $0 \leq i \leq k$. Then for each $0 \leq i \leq k$, the unit root zeta function $L(\phi_i, T)$ is $p$-adic meromorphic.*

Since the unit root pieces $\phi_i$ are not overconvergent any more, Monsky's trace formula cannot be used. Thus, new ideas are needed to prove Dwork's conjecture. By (2.3), the only thing we know is that the following product of these shifted unit root zeta functions

$$(4.4) \qquad L(\phi, T) = \prod_{i=0}^{k} L(\phi_i, \pi^i T)$$

is $p$-adic meromorphic. This implies that the slope zero zeta function $L(\phi_0, T)$ is meromorphic in the non-trivial disk $|T|_\pi < p$. We want to show that it is meromorphic everywhere. Starting with (4.4), we introduced in [19] a limiting approach which can be used to prove the conjecture in the case $U_i$ is of rank one. In the present paper, again starting with (4.4), we introduce an embedding approach to reduce Dwork's conjecture from higher rank case to rank one case.

The most general Hodge-Newton decomposition is due to Katz [10]. We now describe it since it will be used in Katz's isogeny theorem and in our later proof as well.

**Definition 4.11.** *A $\sigma$-module $(M, \phi)$ is called **ordinary** at a break lattice point $(n_1, n_2)$ if $(n_1, n_2)$ is on some basis polygon $P(\vec{e})$ of $\phi$ and if $(n_1, n_2)$ is fibre by fibre a break point of the Newton polygon of $\phi$ at each fibre.*

Note that our condition is stronger than assuming that the Newton polygon goes through the lattice point $(n_1, n_2)$ of the basis polygon $P(\vec{e})$ at each fibre. The lattice point $(n_1, n_2)$ has to be a break point for the Newton polygon fibre by fibre. That is, the Newton polygon has to turn at the lattice point $(n_1, n_2)$. The basis $\vec{e}$ is then called ordinary at the break point $(n_1, n_2)$. For instance, it may happen that the Newton polygon at each fibre only coincides with the basis polygon for a portion (not the whole side) of the smallest slope side. In our terminology, this case is called ordinary at the smallest slope, but not ordinary on the whole smallest slope side. This special case is actually what we will need later on. Theorem 4.9 does not apply to such a "partial" ordinary situation. It is handled by the following more general result.

**Theorem 4.12** (Katz). *Let $(M, \phi)$ be an injective $\sigma$-module of rank $r$. Assume that $\phi$ is ordinary at a break lattice point $(n_1, n_2)$. Then shrinking $X$ if necessary, there is a free $\phi$-stable $A_0$-submodule $M_0$ of rank $n_1$ with free quotient $M/M_0$ such that the Newton slopes of $(M_0, \phi)$ are fibre by fibre the $n_1$ smallest Newton slopes of $(M, \phi)$.*

*Proof.* Again, we just give an outline of the proof. Shrinking $X$ if necessary, we may assume that $\phi$ is fibre by fibre injective.

First, we assume that $n_1 = 1$. In this case, the ordinary condition at $(1, n_2)$ shows that $\phi$ is divisible by $\pi^{n_2}$. Removing the power $\pi^{n_2}$, we may assume that $n_2 = 0$. The desired $A_0$-submodule $M_0$ is (by p156 of [10]) simply the intersection of the images of $\phi^k$:

$$M_0 = \bigcap_{k=1}^{\infty} \phi^k(M^{(\sigma^k)}).$$

The map

$$\phi : M_0^{(\sigma)} \longrightarrow M_0$$

is well defined, $A_0$-linear and invertible. It is shown in [10] (Theorem 2.4.2) that $M_0$ is locally free of rank 1 such that the quotient $M/M_0$ is also locally free. Shrinking $X$ if necessary, we conclude that $(M_0, \phi)$ is a free rank one unit root $\sigma$-module with free quotient $M/M_0$. This proves the theorem for $n_1 = 1$.

To handle the general case $(n_1, n_2)$ with $n_1 > 1$, one applies the above special case to the exterior power $\wedge^{n_1} M$. Shrinking $X$ if necessary, this gives a free rank one $A_0$-submodule $M_0(n_1)$ of $\wedge^{n_1} M$ with free quotient, stable under $\wedge^{n_1} \phi$, such that the Newton slope of $(M_0(n_1), \wedge^{n_1} \phi)$ is fibre by fibre the smallest Newton slope of $\wedge^{n_1} M$. To conclude the proof, one needs to solve the Plücker equation

$$(M_0(n_1), \wedge^{n_1} \phi) = \wedge^{n_1}(M_0, \phi)$$

for some $\phi$-stable rank $n_1$ free $A_0$-submodule $M_0$ of $M$. As shown on p147 in [10], the above Plücker equation can be solved over a suitable extension of $A_0$ and hence can be solved over $A_0$. The slope assertion about $M_0$ follows from the corresponding slope assertion about $M_0(n_1)$. The theorem is proved. □

## 5. An easier case

In this section, we give a complete and self-contained proof of Dwork's conjecture in an easier case, where "easier" means that $(M, \phi)$ is ordinary and the determinant $\sigma$-module of the unit root piece $\phi_i$ of $(M, \phi)$ is overconvergent for all $i$. We use a reduction approach which will prove Dwork's conjecture simultaneously for all allowable unit root $\sigma$-modules, see below.

**Definition 5.1.** Let $C(X)$ be the category of ordinary overconvergent $\sigma$-modules over $X$. Let $C_0(X)$ be the subcategory consisting of those $\sigma$-modules $(M, \phi)$ in $C(X)$ such that each unit root $\sigma$-module $(U_i, \phi_i)$ coming from $(M, \phi)$ has the property that the determinant $\sigma$-module

$$\det(\phi_i) = \wedge^{r_i} \phi_i$$

is overconvergent (namely, isomorphic to an overconvergent one), where $r_i$ denotes the rank of $\phi_i$.

Note that for $\phi \in C_0(X)$, each $\phi_i$ is NOT overconvergent in general. Our condition for $(M, \phi)$ to be in $C_0(X)$ is that the determinant $\sigma$-module $\wedge^{r_i} \phi_i$ is overconvergent for all $i$. This is much weaker than requiring $\phi_i$ itself to be overconvergent if $r_i > 1$. One checks that both $C(X)$ and $C_0(X)$ are closed under direct sum, tensor product, symmetric power and exterior power. Thus, both $C(X)$ and $C_0(X)$ form a tensor category.

**Definition 5.2.** Let $D(X)$ be the set of unit root $\sigma$-modules over $X$ which come from arbitrary slope part in the Hodge-Newton decomposition of elements in $C(X)$.

Let $D_0(X)$ be the set of unit root $\sigma$-modules over $X$ which come from arbitrary slope part in the Hodge-Newton decomposition of elements in $C_0(X)$.

It is clear that we have the inclusions

$$C_0(X) \subset C(X), \quad D_0(X) \subset D(X).$$

We emphasize that the rank one elements in $D_0(X)$ are already overconvergent by our definition. However, a unit root $\sigma$-module in $D_0(X)$ with rank greater than one will NOT be overconvergent in general. We shall see that the set $D_0(X)$ coincides with the tensor category it generates. But the larger set $D(X)$ does not seem to coincide with the tensor category it generates. Dwork's conjecture in the ordinary case says that if $\phi \in D(X)$ (or more generally in the tensor category which $D(X)$ generates), then $L(\phi, T)$ is $p$-adic meromorphic. In this section, we prove the conjecture in the easier case when $\phi$ is in the smaller tensor category $D_0(X)$.

**Theorem 5.3.** *If $\phi \in D_0(X)$, then $L(\phi, T)$ is $p$-adic meromorphic.*

*Proof.* Let $m$ be the unique non-negative integer (or $\infty$ if it does not exist) such that the theorem is true in the open disk $|T|_\pi < p^m$ for all $\phi \in D_0(X)$ but not true in the larger open disk $|T|_\pi < p^{m+1}$ for some $\phi \in D_0(X)$. Clearly, $m \geq 0$. We need to prove that $m = \infty$.

Assume that $m < \infty$. There is then an element $\psi \in D_0(X)$ such that $L(\psi, T)$ is $p$-adic meromorphic in the open disk $|T|_\pi < p^m$ but NOT meromorphic in the larger open disk $|T|_\pi < p^{m+1}$. Namely, $\psi$ is one of the "worst" elements in $D_0(X)$. Let $\psi$ come from the slope $i$ part of some $\phi \in C_0(X)$ for some $i$. That is, $\psi = \phi_i$ for some $0 \leq i \leq k$, where $\phi_0, \phi_1, \cdots, \phi_k$ are the unit root $\sigma$-modules in the Hodge-Newton decomposition of $\phi$. Denote the rank of $\phi_j$ by $r_j$. We may assume that $r_i \geq 1$ as the case $r_i = 0$ has nothing to prove.

Let

$$\rho = \wedge^{r_0}\phi_0 \otimes \wedge^{r_1}\phi_1 \otimes \cdots \otimes \wedge^{r_{i-1}}\phi_{i-1},$$

where the factor $\wedge^{r_j}\phi_j$ is understood to be the identity overconvergent $\sigma$-module $(A, \sigma)$ if $r_j = 0$, by our conventions. This unit root $\sigma$-module $\rho$ is of rank one and overconvergent since $\wedge^{r_j}\phi_j$ is of rank one and overconvergent for all $j$ by our definition of $\phi \in C_0(X)$.

Define a new $\sigma$-module by

$$\varphi = \pi^{-r_1 - \cdots - (i-1)r_{i-1} - i} \wedge^{r_0 + \cdots + r_{i-1} + 1} \phi.$$

This $\sigma$-module is ordinary and overconvergent since $\phi$ is such a $\sigma$-module. Furthermore, $\varphi \in C_0(X)$ since $\phi \in C_0(X)$. Let $\varphi_0, \varphi_1, \cdots$ be the unit root $\sigma$-modules coming from the Hodge-Newton decomposition of $\varphi$. One checks that

$$\varphi_0 = \phi_i \otimes \rho = \psi \otimes \rho.$$

Thus, the slope zero part of $\varphi \otimes \rho^{-1}$ is $\psi$. By (4.4), we have the following decomposition formula:

(5.1) $$L(\varphi \otimes \rho^{-1}, T) = L(\psi, T) \prod_{j \geq 1} L(\varphi_j \otimes \rho^{-1}, \pi^j T).$$

Now, the left side is meromorphic everywhere by Corollary 3.2, because $\varphi \otimes \rho^{-1}$ is overconvergent. On the other hand, for each $j \geq 1$, $\varphi_j \otimes \rho^{-1}$ is in $D_0(X)$ as $\varphi \otimes \rho^{-1} \in C_0(X)$. It follows from our definition of $m$ that the L-function $L(\varphi_j \otimes \rho^{-1}, T)$ is meromorphic in $|T|_\pi < p^m$ for each $j$. Thus, the shifted L-function $L(\varphi_j \otimes \rho^{-1}, \pi^j T)$ is meromorphic in the larger disk $|T|_\pi < p^{m+1}$ for each $j \geq 1$. We conclude that the "only remaining" factor $L(\psi, T)$ in (5.1) is also meromorphic in the larger disk $|T|_\pi < p^{m+1}$. This contradicts our choice of $\psi$. Thus, we must have $m = \infty$. The proof is complete. $\square$

*Remark* 5.4. The above proof shows that each element $\psi$ in $D_0(X)$ can be embedded as the slope zero part of some element $\varphi \otimes \rho^{-1}$ in $C_0(X)$. This is a key property in the above reduction approach. If we could show that each element $\phi$ in $D(X)$ could be embedded as the slope zero part of some element in $C(X)$, then Dwork's conjecture for the more general ordinary case $D(X)$ could be proved in exactly the same way as above. However, the embedding property does not seem to be true for $D(X)$. In next section, we get around the difficulty with the help of our limiting approach. The embedding property for $D_0(X)$ shows that the set $D_0(X)$ is closed under direct sum, tensor product, symmetric power and exterior power. Thus, $D_0(X)$ coincides with the tensor category it generates. The same statement does not seem to be true for the larger set $D(X)$.

## 6. Ordinary case

There are several reasons for the smaller category $D_0(X)$ to be easier. First, it has the embedding property and thus forms a tensor category itself. Second, the rank one elements in $D_0(X)$ are already overconvergent and thus well understood by Corollary 3.2. For the larger set $D(X)$, none of these properties seems to hold. We shall create a larger category containing $D(X)$ which will have the embedding property to some category larger than $C(X)$ but small enough so that our limiting approach can handle its rank one elements. Then the same method as in the easier case goes through with the larger category. This gives us a little more than Dwork conjectured.

**Definition 6.1.** Let $E(X)$ be the set of unit root $\sigma$-modules over $X$ which come from the **slope zero** part in the Hodge-Newton decomposition of elements in $C(X)$. Let $I(X)$ be the set of rank one elements in $E(X)$.

Let $\phi$ and $\psi$ be two elements of $C(X)$. Let $\phi_0$ (resp. $\psi_0$) be the slope zero part of $\phi$ (resp. $\psi$). It is clear that that $\phi_0 \oplus \psi_0$ is the slope zero part of $\phi \oplus \psi$. Similarly, $\phi_0 \otimes \psi_0$ is the slope zero part of $\phi \otimes \psi$. This shows that the set $E(X)$ coincides with the tensor category it generates. But the set $E(X)$ is smaller than the set $D(X)$ which is what we want to study. The subset $I(X)$ is closed under tensor product. If we invert the rank one elements in $E(X)$, then we will get a tensor category containing $D(X)$. In this way, our original ambient category $C(X)$ will have to be enlarged as well but the rank one elements of the enlarged ambient category will be within the range of our limiting approach.

For a rank one element $\phi \in I(X)$, we define the inverse $\phi^{-1}$ (or the contragradient of $\phi$) to be the rank one $\sigma$-module with the same underlying module as $\phi$ but its Frobenius matrix is the inverse (the reciprocal) of the Frobenius matrix of $\phi$. It is clear that $\phi \otimes \phi^{-1}$ is isomorphic to the identity $\sigma$-module $(A_0, \sigma)$. Let $I(X)^{-1}$ be the set consisting of the inverses of elements in $I(X)$. Let $E(X) \otimes I(X)^{-1}$ be the set consisting of the elements $\phi \otimes \psi$, where $\phi \in E(X)$ and $\psi \in I(X)^{-1}$. One checks that the set $E(X) \otimes I(X)^{-1}$ is the same as the tensor category it generates. Furthermore, the set $E(X) \otimes I(X)^{-1}$ consists exactly of the unit root $\sigma$-modules coming from arbitrary slope part in the Hodge-Newton decomposition of elements in $C(X) \otimes I(X)^{-1}$, by the following lemma.

**Lemma 6.2.** *We have the inclusion*
$$E(X) \subset D(X) \subset E(X) \otimes I(X)^{-1}.$$

*Proof.* The first inclusion is obvious. The second inclusion is already implicit in the proof of Theorem 5.3. We repeat the argument here. Let $\psi \in D(X)$. We want to show that $\psi \in E(X) \otimes I(X)^{-1}$.

Let $\psi$ come from the slope $i$ part of some $\phi \in C(X)$. That is, $\psi = \phi_i$ for some $0 \leq i \leq k$, where $\phi_0, \phi_1, \cdots, \phi_k$ are the unit root $\sigma$-modules in the Hodge-Newton decomposition of $\phi$. Denote the rank of $\phi_j$ by $r_j$. We may assume that $r_i \geq 1$ as the case $r_i = 0$ has nothing to prove.

Let
$$\rho = \wedge^{r_0}\phi_0 \otimes \wedge^{r_1}\phi_1 \otimes \cdots \otimes \wedge^{r_{i-1}}\phi_{i-1},$$
where the factor $\wedge^{r_j}\phi_j$ is understood to be the identity $\sigma$-module $\sigma$ if $r_j = 0$, by our convention. This unit root $\sigma$-module $\rho$ is of rank one and it is in $I(X)$ since it is the slope zero part of the element
$$\pi^{-r_1-\cdots-(i-1)r_{i-1}}\wedge^{r_0+\cdots+r_{i-1}}\phi$$
which is in $C(X)$. Define a new $\sigma$-module by
$$\varphi = \pi^{-r_1-\cdots-(i-1)r_{i-1}-i}\wedge^{r_0+\cdots+r_{i-1}+1}\phi.$$

This $\sigma$-module is in $C(X)$ since $\phi \in C(X)$. Let $\varphi_0, \varphi_1, \cdots$ be the unit root $\sigma$-modules coming from the Hodge-Newton decomposition of $\varphi$. One checks that
$$\varphi_0 = \phi_i \otimes \rho = \psi \otimes \rho.$$
This shows that $\psi \otimes \rho$ is in $E(X)$. Since $\rho \in I(X)$, it follows that $\psi \in E(X) \otimes I(X)^{-1}$. The lemma is proved. $\square$

Now, we prove Dwork's conjecture for $\psi \in E(X) \otimes I(X)^{-1}$. Since the new ambient category is $C(X) \otimes I(X)^{-1}$, to apply our embedding approach, we need to assume that the $p$-adic meromorphic continuation holds for $C(X) \otimes I(X)^{-1}$. This is really a slight generalization of the rank one case of Dwork's conjecture called key lemma below.

**Theorem 6.3.** *If $\psi \in E(X) \otimes I(X)^{-1}$, then $L(\psi, T)$ is $p$-adic meromorphic.*

*Proof.* Let $m$ be the unique non-negative integer (or $\infty$ if it does not exist) such that the theorem is true in the open disk $|T|_\pi < p^m$ for all $\psi \in E(X) \otimes I(X)^{-1}$ but not true in the larger open disk $|T|_\pi < p^{m+1}$ for some $\psi \in E(X) \otimes I(X)^{-1}$. Clearly, $m \geq 0$. We need to prove that $m = \infty$.

Assume that $m < \infty$. There is then an element $\psi \in E(X) \otimes I(X)^{-1}$ such that $L(\psi, T)$ is $p$-adic meromorphic in the open disk $|T|_\pi < p^m$ but NOT meromorphic in the larger open disk $|T|_\pi < p^{m+1}$. Namely, $\psi$ is one of the "worst" elements in $E(X) \otimes I(X)^{-1}$. Write
$$\psi = \phi_0 \otimes \rho^{-1},$$
where $\rho \in I(X)$ and $\phi_0$ is the slope zero part of some element $\phi \in C(X)$. Let $\phi_0, \phi_1, \cdots, \phi_k$ denote all the unit root $\sigma$-modules in the Hodge-Newton decomposition of $\phi$. Then, $\psi = \phi_0 \otimes \rho^{-1}, \phi_1 \otimes \rho^{-1}, \cdots, \phi_k \otimes \rho^{-1}$ are all the unit root $\sigma$-modules in the Hodge-Newton decomposition of $\phi \otimes \rho^{-1}$. By (4.4), we have the following decomposition formula:

$$(6.1) \qquad L(\phi \otimes \rho^{-1}, T) = L(\psi, T) \prod_{j \geq 1} L(\phi_j \otimes \rho^{-1}, \pi^j T).$$

For each $j \geq 1$, $\phi_j \otimes \rho^{-1}$ is in $E(X) \otimes I(X)^{-1}$ by Lemma 6.2. It follows from our definition of $m$ that the L-function $L(\phi_j \otimes \rho^{-1}, T)$ is meromorphic in $|T|_\pi < p^m$ for each $j$. Thus, the shifted L-function $L(\phi_j \otimes \rho^{-1}, \pi^j T)$ is meromorphic in the larger disk $|T|_\pi < p^{m+1}$ for each $j \geq 1$. Although the $\sigma$-module $\phi \otimes \rho^{-1}$ is not overconvergent in general, it is the tensor product of the harmless overconvergent $\phi$ with the non-overconvergent but rank one element $\rho^{-1}$ in $I(X)^{-1}$. Our key lemma below shows that its L-function $L(\phi \otimes \rho^{-1}, T)$ is still meromorphic everywhere. We conclude that the "only remaining" factor $L(\psi, T)$ in (6.1) is also meromorphic in the larger disk $|T|_\pi < p^{m+1}$. This contradicts our choice of $\psi$. Thus, we must have $m = \infty$. The proof is complete. □

**Lemma 6.4** (key lemma). *Let $\phi$ be an overconvergent $\sigma$-module over $X$. Let $\rho$ be the slope zero part of some overconvergent $\sigma$-module ordinary at the slope zero side. Assume that $\rho$ has rank one. Then the L-function $L(\phi \otimes \rho^k, T)$ is meromorphic everywhere for every integer $k$. Furthermore, the family of L-functions $L(\phi \otimes \rho^k, T)$ parametrized by integer $k$ in any given residue class modulo $(q-1)$ is a strong family in the sense of [23].*

This result, actually only its first part, is what we need to prove Theorem 6.3 which includes Dwork's conjecture for $D(X)$ as a corollary. It shows that Corollary 3.2 is true for a significantly larger category of $\sigma$-modules. In [19], we proved Lemma 6.4 in the technically simpler setting where Dwork's trace formula applies. In Lemma 9.6 of section 9, we shall further reduce the key lemma from $X$ to the affine $n$-space $\mathbf{A}^n$. But the simpler Dwork trace formula still cannot be used because the Frobenius lifting $\sigma$ is not necessarily the simplest one $X_i \to X_i^q$. One needs to use a sufficiently explicit form of Monsky's trace formula. This is completed in [23].

The above proof is actually constructive. We can use it to obtain an explicit formula for the unit root L-function. For this purpose, we restate a stronger formula version of Lemma 6.4 as follows.

**Lemma 6.5.** *Let $\psi \in C(X) \otimes I(X)^{-1}$. Then $L(\psi, T)$ is $p$-adic meromorphic. Furthermore, there is an effectively determined finite sequence $\Psi_i$ $(0 \leq i \leq m(\psi))$ of nuclear overconvergent $\sigma$-modules depending on $\psi$ such that*

$$L(\psi, T) = \prod_{i=0}^{m(\psi)} L(\Psi_i, T)^{\pm},$$

*where $m(\psi)$ is a positive integer depending on $\psi$.*

This lemma is a consequence of Theorem 7.7 in [23] and (9.8) in section 9 of this paper. The definition of nuclear overconvergent $\sigma$-modules is given in [23]. For our purpose here, one only needs to accept the statement that the L-function of a nuclear overconvergent $\sigma$-module is automatically $p$-adic meromorphic, by the infinite rank version of the Dwork-Monsky trace formula. Lemma 6.5 gives an explicit formula for the unit root L-function in the rank one case in terms of infinite rank nuclear overconvergent $\sigma$-modules. For higher rank case, we also want to get a formula for $L(\phi, T)$ when $\phi \in E(X) \otimes I(X)^{-1}$. First, we need

**Lemma 6.6.** *Let $\phi \in E(X) \otimes I(X)^{-1}$. There is an effectively determined infinite sequence of finite rank $\sigma$-modules:*

$$\psi_i \in C(X) \otimes I(X)^{-1}, \ 0 \leq i < \infty$$

*depending on $\phi$ such that*

$$L(\phi, T) = \prod_{i=0}^{\infty} L(\psi_i, T)^{\pm}, \quad \lim_{i \to \infty} \mathrm{ord}_\pi(\psi_i) = \infty.$$

*Proof.* This follows from Lemma 6.2 and the proof of Theorem 6.3. In fact, write

$$\phi = \varphi_0 \otimes \rho, \ \rho \in I(X)^{-1}, \varphi_0 \in E(X).$$

Let $\varphi_0$ be the slope zero part of some element $\varphi \in C(X)$. Denote the slope $i$ part of $\varphi$ by $\varphi_i$. Then, by the definition of L-functions, we have

(6.2) $$L(\phi, T) = \frac{L(\varphi \otimes \rho, T)}{\prod_{i \geq 1} L(\varphi_i \otimes \rho, \pi^i T)}.$$

The $\sigma$-module $\varphi \otimes \rho$ in the numerator is an element of the better category $C(X) \otimes I(X)^{-1}$ for which Lemma 6.5 applies. Each $\sigma$-module $\varphi_i \otimes \rho$ in the denominator is an element of $E(X) \otimes I(X)^{-1}$, but we do get a non-trivial twisting factor $\pi^i$ ($i \geq 1$). We can apply equation (6.2) to each factor in its denominator. Iterating this procedure, because of the shifting factor $\pi$, we deduce by induction that the lemma is true. □

Combining Lemma 6.5 and Lemma 6.6, we obtain an explicit formula in the higher rank case.

**Theorem 6.7.** *Let $\phi \in E(X) \otimes I(X)^{-1}$. There is an effectively determined infinite sequence of nuclear overconvergent $\sigma$-modules $\Phi_i$ ($0 \leq i < \infty$) depending on $\phi$ such that*

$$L(\phi, T) = \prod_{i=0}^{\infty} L(\Phi_i, T)^{\pm}, \quad \lim_{i \to \infty} \mathrm{ord}_\pi(\Phi_i) = \infty.$$

*Note that each $\Phi_i$ here is in general of infinite rank and the product is in general an infinite product.*

In the more restrictive case of overconvergent F-crystal setting, each $\sigma$-module $\Phi_i$ becomes an overconvergent nuclear F-crystal. The extra integrable connection would allow one to obtain a cohomological formula for each L-function $L(\Phi_i, T)$. Such a cohomological formula explains some cancellation of the zeros and thus leads to improved information about the distribution of the zeros and poles of each $L(\Phi_i, T)$. A more mysterious question is to try to understand if there are any non-trivial differential maps among the different nuclear F-crystals $\Phi_i$ for different $i$. Such non-trivial maps, if they exist, would lead to a different kind of possible cancellation of zeros and poles for the more difficult unit root L-function $L(\phi, T)$.

## 7. Non-ordinary case

In this section, we treat Dwork's conjecture in the more general case where the overconvergent $\sigma$-module $(M, \phi)$ is not necessarily ordinary. The method is similar to previous two sections, except now we need to use Grothendieck's specialization theorem, Katz's isogeny theorem and some additional arguments. Again, we will need to use the key lemma from section 6.

Let $(M, \phi)$ be a $\sigma$-module over $X/\mathbf{F}_q$. Grothendieck's specialization theorem shows that $(M, \phi)$ has a generic Newton polygon. In particular, there are only a finite number of possibilities for the Newton slopes of the $\phi_x$ for all closed points $\bar{x} \in X/\mathbf{F}_q$. By going to a totally ramified finite extension of $R$ if necessary, we may also assume that all Newton slopes of all fibres $\phi_x$ are non-negative integers. For an integer $i \geq 0$, recall that we defined the slope $i$ part of the characteristic polynomial of $(M, \phi)$ at $x$ to be the product

$$\det_i(\phi_x, T) = \prod_{\mathrm{ord}_{\pi^{\deg(x)}}(\alpha) = i} (1 - \alpha T),$$

where $\alpha$ runs over the reciprocal roots of slope $i$ of the characteristic polynomial $\det(I - \phi_x^{\deg(x)} T)$ defined with respect to the valuation $\mathrm{ord}_{\pi^{\deg(x)}}$. It satisfies the product relation

(7.1) $$\det(I - \phi_x^{\deg(x)} T) = \prod_{i \geq 0} \det_i(\phi_x, T).$$

Define

$$L_i(\phi, T) = \prod_{\bar{x} \in X/\mathbf{F}_q} \frac{1}{\det_i(\phi_x, T^{\deg(x)})} \in R[[\pi^i T]],$$

where $\bar{x}$ runs over the closed points of $X/\mathbf{F}_q$. This is called the slope $i$ zeta function of $(M, \phi)$. It is trivially analytic in $|T|_\pi < p^i$. By (7.1), we get the decomposition formula:

$$L(\phi, T) = \prod_{i \geq 0} L_i(\phi, T).$$

Conjecture 2.4 in this case becomes

**Conjecture 7.1** (Dwork). *Let $(M, \phi)$ be an overconvergent $\sigma$-module over $X/\mathbf{F}_q$ with integral Newton slopes. Then for each integer $i \geq 0$, the slope $i$ zeta function $L_i(\phi, T)$ is p-adic meromorphic.*

To adapt our proof in the ordinary case to the more general non-ordinary case, we need to use Katz's isogeny theorem [10] for an injective $\sigma$-module and to understand the stability of the overconvergence under certain weak Newton-Hodge isogenies. Recall that $\phi$ is injective if and only if the determinant $\det \phi$ is non-zero.

**Theorem 7.2** (Katz). *Let $(M, \phi)$ be an injective $\sigma$-module over $X$ of rank $r$ with integral Newton slopes. Then, shrinking $X$ if necessary, the $\sigma$-module $(M, \phi)$ is isogenous to an **ordinary** $\sigma$-module $(M', \phi')$. That is, the $\sigma$-module $(M', \phi')$ admits a $\phi'$-stable filtration of free $A_0$-modules with free quotients:*

$$0 \subset M_0' \subset M_1' \subset M_2' \subset \cdots$$

*such that the quotient $(M_i'/M_{i-1}', \phi')$ is a $\sigma$-module of the form $(U_i, \pi^i \phi_i')$, where $\phi_i'$ is a unit root $\sigma$-module for each $i \geq 0$.*

*Proof.* This theorem follows from the remark on page 154 of [10]. We give a sketch of the proof as we will need another result derived from the proof. By Grothendieck's specialization theorem, we can shrink $X$ if necessary and assume that the generic Newton polygon of $(M, \phi)$ coincides with the Newton polygon of each fibre $\phi_x$. Namely, the Newton polygon of $\phi_x$ is a constant Newton polygon as $\bar{x}$ moves on the variety $X/\mathbf{F}_q$. In particular, $\phi$ is injective fibre by fibre. Let $s$ be the smallest non-trivial Newton slope of $\phi$, which is a non-negative integer by our assumption. Lemma 4.6 implies that for all integers $k \geq 0$, we have
$$\phi^{k+r-1} \equiv 0 \pmod{\pi^{ks}}.$$
Let $\nu = (r-1)\text{ord}_\pi \det \phi$, which is a finite integer since $\phi$ is injective. Now, the inverse of the matrix for $\phi^{r-1}$ has entries in $\pi^{-\nu} A_0$. It follows from the above congruence that for all integers $k \geq 0$, we have
$$\phi^k \equiv 0 \pmod{\pi^{ks-\nu}}.$$
Define an $A_0$-module $M'$ with
$$M \subset M' \subset \pi^{-\nu} M$$
by
$$M' = \sum_{k \geq 0} \frac{\phi^k}{\pi^{ks}}(M^{(\sigma^k)}),$$
Note that the summand for $k=0$ gives $M$. It is clear that $\phi$ maps $M'$ into itself and is $\sigma$-linear. Furthermore, $(M', \phi)$ is divisible by $\pi^s$. That is, $\phi/\pi^s$ also maps $M'$ into itself. Denote the pair $(M', \phi)$ by $(M', \phi')$.

We now prove that $M'$ is a free $A_0$-module at the expense of further shrinking $X$. Let $L$ be the localization of $A_0$ at the prime ideal $\pi A_0$. Namely, we invert all elements in $A_0$ which are not divisible by $\pi$. One checks that $L$ is a discrete valuation ring with all ideals being of the form $(\pi^i)$. Since $L$ is a principal ideal domain, the inclusion
$$M \otimes L \subset M' \otimes L \subset \pi^{-\nu}(M \otimes L)$$
shows that $M' \otimes L$ is a free $L$-module of rank $r$. Furthermore, $(M \otimes L, \phi)$ is "isogenous" to $(M' \otimes L, \phi')$ which is a $\sigma$-module over $L$ divisible by $\pi^s$. That is, $\phi'$ also maps $M' \otimes L$ into itself and is $\sigma$-linear. Let $\vec{e}$ be a fixed basis of the free $A_0$-module $M$. By scalar extension, $\vec{e}$ is a basis of the free $L$-module $M \otimes L$. Similarly, $\pi^{-\nu}\vec{e}$ is a basis of the free $A_0$-module $\pi^{-\nu}M$ (and also the free $L$-module $\pi^{-\nu}M \otimes L$). Let $\vec{f}$ be a fixed basis of the free $L$-module $M' \otimes L$. Let $E$ be the transition matrix from $\vec{e}$ to $\vec{f}$. That is,
$$\vec{f} = \vec{e}E,$$
where the matrix $E$ has entries in $L$ with non-zero determinant. We deduce that there is an element $g \in A_0$ not divisible by $\pi$ such that both $gE$ and $gE^{-1}$ have entries in $A_0 \otimes K$, where $K$ is the quotient field of the constant ring $R$. This means that if we remove the hypersurface section $g = 0$ from $X$ if necessary, then on the smaller smooth affine variety $X - \{g = 0\}$, the $\sigma$-module $(M, \phi)$ is isogenous to a $\sigma$-module $(M', \phi')$ which is divisible by $\pi^s$. The argument so far works for overconvergent $\sigma$-modules as well. That is, if $(M, \phi)$ is overconvergent, then we can take $(M', \phi')$ to be overconvergent. One simply replaces $A_0$ in the above proof by $A$.

We now finish the proof by induction on the number of distinct Newton slopes, following p157 in [10]. This part needs the Hodge-Newton decomposition. The inclusion of $(M, \phi)$ into $(M', \phi')$ is an isogeny such that $(M', \phi')$ is divisible by $\pi^s$. Since $(M, \phi)$ and $(M', \phi')$ have the same Newton polygon, it follows that $(M', \phi')$ is ordinary at the smallest slope $s$. If $(M, \phi)$ has only one Newton slope, then $(M', \phi')$ is ordinary everywhere and we are already done. In general, $(M', \phi')$ will be ordinary at some lattice break point on the smallest slope $s$ side. We can apply Katz's Hodge-Newton decomposition (Theorem 4.12) to this break lattice point. This produces a short exact sequence

$$0 \longrightarrow (M'_s, \phi') \longrightarrow (M', \phi') \longrightarrow (M'/M'_s, \phi') \longrightarrow 0,$$

where $M'_s$ denotes the pure slope $s$ part of $M'$ and the Newton slopes of the quotient $M'/M'_s$ are strictly larger than $s$. By the induction hypothesis applied to $(M'/M'_s, \phi')$, we get an isogeny

$$((M'/M'_s)'', \phi'') \longrightarrow (M'/M'_s, \phi')$$

whose source satisfies all the conclusions of the theorem. Taking the pull back by this map of the above extension of $(M'/M'_s, \phi')$ by $(M'_s, \phi')$, we get an extension

$$0 \longrightarrow (M'_s, \phi') \longrightarrow ? \longrightarrow ((M'/M'_s)'', \phi'') \longrightarrow 0.$$

This exact sequence and the inductively given filtration on $((M'/M'_s)'', \phi'')$ provide the desired filtration of the theorem. The proof is complete. □

The above formulation of Katz's isogeny theorem is best possible. However, in the case of curves and locally free version, Katz [10] showed that no shrinking is necessary if in addition $\phi$ has a constant Newton polygon fibre by fibre. Katz asked if the isogeny theorem for locally free version is still true without shrinking $X$ for higher dimensional $X$ assuming that $\phi$ has a constant Newton polygon fibre by fibre. This question seems still open.

Combining Theorem 7.2 and Lemma 4.5, we deduce

**Corollary 7.3.** *Let $(M, \phi)$ be a $\sigma$-module over $X$ with integral Newton slopes. Then, shrinking $X$ if necessary, the $\sigma$-module $(M, \phi)$ is isogenous to a $\sigma$-module $(M', \phi')$, where $(M', \phi')$ admits a $\phi'$-stable filtration of free $A_0$-modules with free quotients:*

$$0 \subset N_0 \subset M'_0 \subset M'_1 \subset M'_2 \subset \cdots$$

*such that the restriction of $\phi'$ to $N_0$ is nilpotent and the quotient $(M'_i/M'_{i-1}, \phi')$ is a $\sigma$-module of the form $(U_i, \pi^i \phi'_i)$, where $\phi'_i$ is a unit root $\sigma$-module for each $i \geq 0$.*

Note that Theorem 7.2 and Corollary 7.3 are false over the smaller overconvergent ring $A$ because the Hodge-Newton decomposition will destroy the overconvergent property. However, we do have the following weaker result about stability of overconvergent $\sigma$-modules under an isogeny map. This result is crucial for the harder non-ordinary version of the key lemma. It holds only for the smallest slope part because we cannot apply the Hodge-Newton decomposition in order to keep the overconvergence property.

**Lemma 7.4.** *Let $(M, \phi)$ be an injective overconvergent $\sigma$-module over $X/\mathbf{F}_q$ with integral Newton slopes. Then, shrinking $X$ if necessary, the $\sigma$-module $(M, \phi)$ is isogenous to an **overconvergent** $\sigma$-module $(M', \phi')$ which is ordinary for the whole side of the smallest (non-trivial) Newton slope.*

*Proof.* As observed in the proof of Theorem 7.2, if we do not use the Hodge-Newton decomposition, we deduce the weaker result that shrinking $X$ if necessary, the overconvergent $\phi$ is isogenous to an overconvergent $\phi'$ which is ordinary at the smallest slope, although not necessarily ordinary on the whole side of the smallest slope. Removing a power of $\pi$, we can assume that $\phi$ is overconvergent and ordinary at slope zero. Furthermore, $\phi$ has a constant Newton polygon. We want to show that shrinking $X$ if necessary, $\phi$ is isogenous to an overconvergent $\sigma$-module which is ordinary for the whole horizontal side. Let $h_0 \geq 1$ be the length of the slope zero side of the Newton polygon of $\phi$. Let $B$ be the matrix of $\phi$ with respect to some basis. It has entries in $A$. By Theorem 7.2, shrinking $X$ if necessary, there is an invertible matrix $E$ with entries in $A_0 \otimes \mathbf{Q}$ such that $E^{-1}BE^\sigma$ has the form

$$(7.2) \qquad E^{-1}BE^\sigma = \begin{pmatrix} B_{00} & B_{01} \\ 0 & \pi B_{11} \end{pmatrix}$$

where $B_{00}$ is an invertible $h_0 \times h_0$ matrix with entries in $A_0$ and all $B_{ij}$ have entries in $A_0$. Note that the matrix in (7.2) is not overconvergent. We now deform it to make it overconvergent.

Choose a sufficiently large integer $k$ such that

$$\pi^k E \equiv 0 \pmod{\pi}, \quad \pi^k E^{-1} \equiv 0 \pmod{\pi}.$$

That is, the denominators of $E$ and $E^{-1}$ are at most $\pi^{k-1}$. Choose an overconvergent invertible matrix $E_1$ with entries in $A \otimes \mathbf{Q}$ such that

$$E_1 \equiv E \pmod{\pi^k}, \quad E_1^{-1} \equiv E^{-1} \pmod{\pi^k},$$

where the congruence is taken in $A_0$. Then, the new matrix $C = E_1^{-1}BE_1^\sigma$ is overconvergent and it has the block form

$$(7.3) \qquad C = E_1^{-1}BE_1^\sigma = \begin{pmatrix} C_{00} & C_{01} \\ \pi C_{10} & \pi C_{11} \end{pmatrix},$$

where $C_{00}$ is an invertible $h_0 \times h_0$ matrix with entries in $A$ and all $C_{ij}$ have entries in $A$. The overconvergent $\sigma$-module $C$ defined by the matrix $C$ in (7.3) is isogenous to the overconvergent $\sigma$-module $\phi$.

The $\sigma$-module $C$ in (7.3) is not yet ordinary on the whole horizontal side. To finish the proof, it suffices to prove the claim that (shrinking $X$ if necessary) there is a basis $\vec{e}$ of the $\sigma$-module $C$ over $A \times \mathbf{Q}$ such that the new matrix $D$ of $C$ under $\vec{e}$ is of the form

$$(7.4) \qquad D = \begin{pmatrix} D_{00} & \pi D_{01} \\ \pi D_{10} & \pi D_{11} \end{pmatrix},$$

where $D_{00}$ is an invertible $h_0 \times h_0$ matrix over $A$ and all $D_{ij}$ have entries in $A$. To prove the claim, we use induction. We may assume that the rank $r$ of $M$ is greater than $h_0$ (otherwise, there is nothing to prove). Write the matrix $C$ in column form

$$C = (C_1, \cdots, C_r).$$

Since $\det C_{00}$ is invertible in $A$, the block form in (7.3) shows that shrinking $X$ if necessary, we can write the last column $C_r$ of $C$ as an $A$-linear combination of the first $h_0$ columns modulo $\pi$. That is, there are elements $a_i \in A$ ($1 \leq i \leq h_0$) such that

$$(7.5) \qquad C_r \equiv a_1 C_1 + \cdots + a_{h_0} C_{h_0} \pmod{\pi}.$$

Since

$$C^\sigma = C^{-1}CC^\sigma,$$

we deduce that the $\sigma$-module $C$ is isogenous to the $\sigma$-module $C^\sigma$ which is still overconvergent and of the block form (7.3). Thus, we may replace $C$ by $C^\sigma$ in our proof. Let $E_2$ be the following overconvergent elementary matrix

$$E_2 = \begin{pmatrix} I_{h_0} & E_3 \\ 0 & I_{r-h_0} \end{pmatrix},$$

where $E_3$ is the $h_0 \times (r - h_0)$ matrix all of whose columns are zero except for the last column which is the transpose of the row vector $(-a_1, -a_2, \cdots, -a_{h_0})$. The product matrix $CE_2$ is then the matrix obtained from $C$ by performing the following elementary column operations: multiplying the $i$-th column $(1 \leq i \leq h_0)$ of $C$ by $-a_i$ and adding it to the last column of $C$. Equation (7.5) shows that the last column of $CE_2$ is divisible by $\pi$. Thus, the last column of the new matrix

$$E_2^{-1} C^\sigma E_2^\sigma = E_2^{-1}(CE_2)^\sigma$$

is also divisible by $\pi$. Furthermore, the block form of $E_2$ shows that $E_2^{-1} C^\sigma E_2^\sigma$ has the same block form as (7.3). It is still overconvergent and its last column is now divisible by $\pi$. Applying the above procedure to the columns indexed by $\{h_0 + 1, \cdots, r - 1\}$, we conclude that the claim is true. The lemma is proved. $\square$

With these preparations, we can now return to the proof of Dwork's conjecture in non-ordinary case. We say that a $\sigma$-module $(M, \phi)$ is isogenous to an overconvergent $\sigma$-module $(M', \phi')$ if $(M, \phi)$ is isogenous to the $\sigma$-module obtained from the overconvergent $\sigma$-module $(M', \phi')$ by base extension $A \to A_0$.

**Definition 7.5.** Let $C_1(X)$ be the category of **ordinary** (not necessarily overconvergent) $\sigma$-modules over $X/\mathbf{F}_q$ which are isogenous (without shrinking $X$) to **overconvergent** (not necessarily ordinary) $\sigma$-modules over $X/\mathbf{F}_q$. Let $D_1(X)$ be the set of unit root $\sigma$-modules over $X$ coming from arbitrary slope part in the Hodge-Newton decomposition of elements in $C_1(X)$. Let $E_1(X)$ be the set of unit root $\sigma$-modules over $X$ coming from the smallest non-trivial slope part in the Hodge-Newton decomposition of elements in $C_1(X)$. Let $I_1(X)$ be the set of rank one elements in $E_1(X)$.

It is clear that both $E_1(X)$ and $I_1(X)$ are closed under tensor product. One checks that the set $E_1(X) \otimes I_1(X)^{-1}$ is the same as the tensor category it generates. By our definition, we have the inclusion relations:

$$C(X) \subset C_1(X), \ D(X) \subset D_1(X), \ E(X) \subset E_1(X), \ I(X) \subset I_1(X).$$

The following result can be proved in the same way as was Lemma 6.2.

**Lemma 7.6.** *We have the inclusion*

$$E_1(X) \subset D_1(X) \subset E_1(X) \otimes I_1(X)^{-1}.$$

The next result is an extension of our key lemma (Lemma 6.4) from ordinary case to non-ordinary case. It follows directly from Lemma 6.4, Lemma 4.5 and Lemma 7.4.

**Lemma 7.7** (non-ordinary key lemma)**.** *If $\eta$ is an overconvergent $\sigma$-module over $X$ and if $\rho \in I_1(X)$, then the L-function $L(\eta \otimes \rho^k, T)$ is meromorphic everywhere for every integer $k$. Furthermore, the family of L-functions $L(\eta \otimes \rho^k, T)$ parametrized by integer $k$ in any given residue class modulo $(q-1)$ is a strong family in the sense of* [23].

Note that Lemma 7.7 is also true if we replace $\eta$ by an element in $C_1(X)$. This is because $\rho$ is of rank one and an isogeny map does not change the L-function. Using this fact and Lemma 7.7, we deduce, as in Theorem 6.2, the following result.

**Theorem 7.8.** *If $\psi \in E_1(X) \otimes I_1(X)^{-1}$, then $L(\psi, T)$ is p-adic meromorphic.*

This result easily implies Dwork's conjecture in the non-ordinary case.

**Theorem 7.9.** *Let $(M, \phi)$ be an overconvergent $\sigma$-module over $X$ with integral Newton slopes. Then for each integer $i \geq 0$, the slope $i$ zeta function $L_i(\phi, T)$ is p-adic meromorphic.*

*Proof.* We use induction on the dimension of $X$. The theorem is trivially true if $\dim(X) = 0$. Assume that $\dim(X) > 0$. Assume that the theorem is true for all overconvergent $\sigma$-modules $(M, \phi)$ over all smooth affine varieties of dimension smaller than $\dim(X)$. By Katz's isogeny theorem and Lemma 4.5, we may assume that there is a smooth affine Zariski dense subset $Y$ of $X$ such that $\phi$ restricted to $Y$ is isogenous to an element in $C_1(Y)$. Let $L_i'(\phi, T)$ (resp. $L_i''(\phi, T)$) be the slope $i$ zeta function of $\phi$ restricted to $Y$ (resp. restricted to $X - Y$). By induction, one easily shows that $L_i''(\phi, T)$ is meromorphic. By Lemma 7.6 and Theorem 7.8, $L_i'(\phi, T)$ is meromorphic. It follows that the total slope $i$ zeta function
$$L_i(\phi, T) = L_i'(\phi, T) L_i''(\phi, T)$$
is meromorphic. The proof is complete.   $\square$

Similar to Theorem 6.7, we can derive an explicit formula for the pure slope L-function $L_i(\phi, T)$ in terms of nuclear overconvergent $\sigma$-modules. We have

**Theorem 7.10.** *Let $(M, \phi)$ be an overconvergent $\sigma$-module over $X$ with integral Newton slopes. For each integer $i \geq 0$, there is an effectively determined infinite sequence of "constructible" nuclear overconvergent $\sigma$-modules $\Phi_{i,j}$ $(0 \leq j < \infty)$ depending on $\phi$ such that*
$$L_i(\phi, T) = \prod_{j=0}^{\infty} L(\Phi_{i,j}, T)^{\pm}, \quad \lim_{j \to \infty} \mathrm{ord}_\pi(\Phi_{i,j}) = \infty.$$

In this theorem, "constructible" means that the $\sigma$-module $\Phi_{i,j}$ is a direct sum of finitely many nuclear overconvergent $\sigma$-modules over various subvarieties of $X$. This direct sum is caused by the combinatorial fact that we have to shrink the variety $X$ at various points in the proof.

## 8. The general case

In this section, we treat Dwork's conjecture in full generality by working with the pure slope zeta function attached to a higher power (iterate) of an overconvergent $\sigma$-module over $X/\mathbf{F}_q$.

Let $k$ be a positive integer and let $(M, \phi)$ be a $\sigma$-module. The $k$-th iterate $(M, \phi^k)$ is a $\sigma^k$-module over $X/\mathbf{F}_{q^k}$ but not a $\sigma$-module over $X/\mathbf{F}_q$. Nevertheless, we can still define its L-function over $X/\mathbf{F}_q$, called the $k$-th power L-function of $\phi$, to be the following Euler product
$$L(\phi^k, T) = \prod_{\bar{x} \in X/\mathbf{F}_q} \frac{1}{\det(I - \phi_x^{k\deg(x)} T^{\deg(x)})},$$

where $\bar{x}$ runs over the closed points of $X/\mathbf{F}_q$ and $x$ is the Teichmüller lifting of $\bar{x}$ defined with respect to $\sigma$. Namely, we are raising each characteristic root of the Euler factors of $L(\phi, T)$ to its $k$-th power. We emphasize that $L(\phi^k, T)$ is different from the L-function of $\phi^k$ viewed as a $\sigma^k$-module over $X/\mathbf{F}_{q^k}$ since the degree of a point is different when viewed over $\mathbf{F}_q$ and over its $k$-th extension $\mathbf{F}_{q^k}$. More generally, if we are given two $\sigma$-modules $\phi_1$ and $\phi_2$ over $X/\mathbf{F}_q$, two non-negative integers $k_1$ and $k_2$, then we can define the L-function of the "tensor product" $\phi_1^{k_1} \otimes \phi_2^{k_2}$ of their iterates as follows

$$L(\phi_1^{k_1} \otimes \phi_2^{k_2}, T) = \prod_{\bar{x} \in X/\mathbf{F}_q} \frac{1}{\det(I - \phi_{1,x}^{k_1 \deg(x)} \otimes \phi_{2,x}^{k_2 \deg(x)} T^{\deg(x)})},$$

where $\phi_{1,x}$ denotes the fibre of $\phi_1$ at $x$ and $\phi_{2,x}$ denotes the fibre of $\phi_2$ at $x$.

Again, by Grothendieck's specialization theorem, we can shrink $X$ and go to a totally ramified finite extension of $R$, if necessary. In this way, we may assume that the Newton slopes of $\phi$ are fibre by fibre non-negative integers. For an integer $i \geq 0$, we define

$$\det_i(\phi_x^k, T) = \prod_{\text{ord}_\pi(\alpha) = i \deg(x)} (1 - \alpha^k T),$$

where $\alpha$ runs over the reciprocal roots of slope $i$ of the characteristic polynomial $\det(I - \phi_x^{\deg(x)} T)$ defined with respect to the valuation $\text{ord}_{\pi^{\deg(x)}}$. Define

$$L_i(\phi^k, T) = \prod_{\bar{x} \in X/\mathbf{F}_q} \frac{1}{\det_i(\phi_x^k, T^{\deg(x)})} \in R[[\pi^{ik} T]],$$

where $\bar{x}$ runs over the closed points of $X/\mathbf{F}_q$. This is called the "slope $i$ and $k$-th power" zeta function of $(M, \phi)$ over $X/\mathbf{F}_q$. It is trivially analytic in $|T|_\pi < p^{ik}$. We have the decomposition formula:

$$L(\phi^k, T) = \prod_{i \geq 0} L_i(\phi^k, T).$$

This is a finite product by Grothendieck's specialization theorem. Although the $k$-th iterate $\phi^k$ does not have a $\sigma$-module structure over $X/\mathbf{F}_q$, it can be viewed as a virtual $\sigma$-module over $X/\mathbf{F}_q$ in the sense that it can be written as an integral combination of $\sigma$-modules over $X/\mathbf{F}_q$ with positive and negative coefficients. Namely, we have the following decomposition formula [19]:

(8.1) $$L(\phi^k, T) = \prod_{j \geq 1} L(\text{Sym}^{k-j} \phi \otimes \wedge^j \phi, T)^{(-1)^{j-1} j}.$$

More generally, by [19], we have

(8.2) $$L(\phi_1^{k_1} \otimes \phi_2^{k_2}, T) = \prod_{j_1, j_2 \geq 1} L(\text{Sym}^{k_1 - j_1} \phi_1 \otimes \wedge^{j_1} \phi_1 \otimes \text{Sym}^{k_2 - j_2} \phi_2 \otimes \wedge^{j_2} \phi_2, T)^{(-1)^{j_1 + j_2 - 2} j_1 j_2}.$$

These decomposition formulas actually hold fibre by fibre for each Euler factor. Taking the slope $i$ (or $ik$) part at each Euler factor, we deduce from (8.1) that

(8.3) $$L_i(\phi^k, T) = \prod_{j \geq 1} L_{ik}(\text{Sym}^{k-j} \phi \otimes \wedge^j \phi, T)^{(-1)^{j-1} j}.$$

Applying Theorem 7.9 to the right side of (8.3), we deduce the following theorem which is the full generality as conjectured by Dwork.

**Theorem 8.1.** *Let $(M, \phi)$ be an overconvergent $\sigma$-module over $X/\mathbf{F}_q$ with integral Newton slopes. Then for each integer $i \geq 0$ and each integer $k \geq 1$, the slope $i$ and $k$-th power L-function $L_i(\phi^k, T)$ is p-adic meromorphic.*

Note that the $k$-th power zeta functions $L(\phi^k, T)$ and $L_i(\phi^k, T)$ can also be defined for negative integer $k$ with a straightforward modification of our definition. One only needs to remove those zero characteristic roots $\alpha$ in the Euler product definition since $\alpha^k$ is undefined when $\alpha = 0$ and $k$ is negative. With this generalized definition of $L(\phi^k, T)$ and $L_i(\phi^k, T)$, Theorem 8.1 holds for negative integers as well and hence holds for all integers $k$.

Similarly, Theorem 7.10 carries over to the $k$-th power situation. We have,

**Theorem 8.2.** *Let $(M, \phi)$ be an overconvergent $\sigma$-module over $X$ with integral Newton slopes. Then, there is an effectively determined infinite sequence of "constructible" nuclear overconvergent $\sigma$-modules $\Phi_{i,j,k}$ ($0 \leq i, j, k < \infty$) depending on $\phi$ such that for all non-negative integers $i$ and $k$, we have*

$$L_i(\phi^k, T) = \prod_{j=0}^{\infty} L(\Phi_{i,j,k}, T)^{\pm}, \quad \lim_{j \to \infty} \mathrm{ord}_\pi(\Phi_{i,j,k}) = \infty.$$

To conclude this section, we raise a question on the possible uniform variation of unit root L-functions. Suppose that $\phi$ is an overconvergent $\sigma$-module which is ordinary at slope zero. Let $\phi_0$ be the slope zero part of $\phi$. Denote the rank of $\phi_0$ by $r_0$. We now know that for each positive integer $k$, the unit root L-function $L(\phi_0^k, T)$ is $p$-adic meromorphic. In the case that the rank $r_0$ is one, a stronger result is known. In fact, the key lemma shows that if $r_0 = 1$, then $L(\phi_0^k, T)$ forms a strong family when $k$ varies in a fixed residue class modulo $(q-1)$. That is, the family is the quotient of two uniformly continuous and uniformly entire functions for $k$ in each residue class modulo $(q-1)$. We want to know if the same statement is true in higher rank case $r_0 > 1$.

First, we do have a similar uniform continuity result. Let $N(q, r_0)$ be the order of the finite group $\mathrm{GL}_{r_0}(\mathbf{F}_q)$. That is,

$$N(q, r_0) = (q^{r_0} - 1)(q^{r_0} - q) \cdots (q^{r_0} - q^{r_0 - 1}).$$

For any element $A \in \mathrm{GL}_{r_0}(R)$, it is easy to see that

(8.4) $$A^{N(q, r_0)} \equiv I_{r_0} \pmod{\pi},$$

where $I_{r_0}$ denotes the $r_0 \times r_0$ identity matrix. The congruence in (8.4) implies that whenever $k_1$ and $k_2$ are two integers satisfying

(8.5) $$k_1 \equiv k_2 \pmod{N(q, r_0) p^m},$$

we have

(8.6) $$A^{k_1} \equiv A^{k_2} \pmod{\pi^{m+1}}.$$

This congruence can be improved somewhat if we take into the consideration of the ramification degree of $R$. We shall however not do so for simplicity. From the Euler product definition of the L-function, we deduce from (8.6) that if $k_1$ and $k_2$ are two positive integers satisfying (8.5), then

(8.7) $$L(\phi_0^{k_1}, T) \equiv L(\phi_0^{k_2}, T) \pmod{\pi^{m+1}}.$$

Thus, we do have the uniform continuity result for the family of unit root L-functions $L(\phi_0^k, T)$ parametrized by positive integer $k$ if we restrict $k$ to vary in a fixed residue class modulo $N(q, r_0)$.

For each non-negative rational number $s$, let $D_s(\phi_0, k)$ denote the total degree of the slope $s$ part of the meromorphic function $L(\phi_0^k, T)$, where the total (sum) degree means the sum of the degree of the numerator and the degree of the denominator.

**Question 8.3.** *For each fixed slope $s \geq 0$, can the total degree $D_s(\phi_0, k)$ be uniformly bounded for all integers $k$?*

The key lemma shows that the answer is positive if the rank $r_0 = 1$. A positive answer for $r_0 > 1$ would imply that a weak version of the Gouvêa-Mazur type conjecture extends to $L(\phi_0^k, T)$ for higher rank $\phi_0$. Instead of using the total degree, one could use the usual (difference) degree of a rational function. For each non-negative rational number $s$, let $d_s(\phi_0, k)$ denote the (difference) degree of the slope $s$ part of the meromorphic function $L(\phi_0^k, T)$. This can be a negative integer. Similarly, we can ask

**Question 8.4.** *For each fixed slope $s \geq 0$, can the (difference) degree $d_s(\phi_0, k)$ be uniformly bounded for all integers $k$?*

Since
$$|d_s(\phi_0, k)| \leq D_s(\phi_0, k),$$
a positive answer to Question 8.3 implies a positive answer to Question 8.4. Thus, Question 8.4 is weaker. Again the answer to Question 8.4 is positive for $r_0 = 1$ by the key lemma. We know a tiny bit more in the case of Question 8.4. For $s = 0$, it is easy to see that the (difference) degree function $d_0(\phi_0, k)$ is always bounded. In fact, $d_0(\phi_0, k)$ is just the (difference) degree of the rational function $L(\phi_0^k, T)$ modulo $\pi$. A uniform bound for $d_0(\phi_0, k)$ then follows from the continuity congruence in (8.7). However, the author does not know if the total degree function $D_0(\phi_0, k)$ (for slope $s = 0$) is bounded if $r_0 > 1$, even if $\phi_0$ itself is overconvergent.

If one replaces the $k$-th iterate $\phi_0$ by the $k$-th tensor power $\phi_0^{\otimes k}$ or by the $k$-th symmetric power $\mathrm{Sym}^k \phi_0$, then the corresponding uniform bound questions have a negative answer if $r_0 > 1$. The simplest counter-example is to take $\phi_0$ to be the direct sum of $r_0$ copies of the identity overconvergent $\sigma$-module $(A, \sigma)$. For the $k$-th iterate $\phi_0^k$, although the author does not have a counter-example at this point, he is inclined to believe that Question 8.3 would in general have a negative answer for $r_0 > 1$. Of course, he would be surprised but very pleased if the answer turned out to be positive. For Question 8.4, there is a slightly larger chance for a positive answer, but it is still doubtful at least for $r_0 > 2$.

## 9. REDUCTION TO THE BASE SCHEME $\mathbf{A}^n$

In this section, we further reduce the key lemma from a smooth affine base scheme $X/\mathbf{F}_q$ to the simplest base scheme $\mathbf{A}^n/\mathbf{F}_q$. There are several ways to achieve this reduction. Here we use Dwork's fundamental F-crystal together with some of the ideas from our previous limiting approach.

First, we lift the given overconvergent $\sigma$-module $(M, \phi)$ over $X/\mathbf{F}_q$ to an overconvergent $\sigma$-module over $\mathbf{A}^n/\mathbf{F}_q$. For this purpose, we need to lift $\sigma$ and $A$. Namely, we first lift the identity overconvergent $\sigma$-module $(A, \sigma)$ from $X/\mathbf{F}_q$ to $\mathbf{A}^n/\mathbf{F}_q$. Let $R[X]^\dagger$ be the overconvergent power series ring $R[X_1, \cdots, X_n]^\dagger$. The context should distinguish the two uses of the same letter $X$, one for the variety $X$ and the other for the variable $X = (X_1, \cdots, X_n)$. Recall from the construction in section 2, we have an exact sequence of $R$-algebras:

$$0 \longrightarrow I \longrightarrow R[X]^\dagger \xrightarrow{\mu} A \longrightarrow 0,$$

where $\mu$ is the reduction map of $R[X]^\dagger$ modulo the ideal $I$ of $R[X]^\dagger$. One can lift (non-uniquely) $\sigma$ from $A$ to $R[X]^\dagger$. To do so, it suffices to choose elements $\sigma(X_i)$ of $R[X]^\dagger$ satisfying the following two conditions:

$$\sigma(X_i) \equiv X_i^q \pmod{\pi}, \quad \sigma(X_i) \in \mu^{-1} \circ \sigma \circ \mu(X_i).$$

Such a choice of elements is clearly possible since

$$\sigma \circ \mu(X_i) \equiv X_i^q \pmod{(\pi, I)}.$$

These elements $\sigma(X_i)$ extend to an $R$-algebra endomorphism (still denoted by $\sigma$) of $R[X]^\dagger$ such that the diagram

$$\begin{array}{ccccc} R[X]^\dagger & \xrightarrow{\mu} & A & \longrightarrow & 0 \\ \downarrow \sigma & & \downarrow \sigma & & \\ R[X]^\dagger & \xrightarrow{\mu} & A & \longrightarrow & 0 \end{array}$$

commutes. That is, we have

$$\sigma \circ \mu = \mu \circ \sigma.$$

Iterating this equation, we deduce that for each positive integer $d$,

$$\sigma^d \circ \mu = \mu \circ \sigma^d.$$

In particular, we have

$$\sigma(I) \subset I.$$

The pair $(R[X]^\dagger, \sigma)$ is our overconvergent lifting of $(A, \sigma)$ from $X/\mathbf{F}_q$ to $\mathbf{A}^n/\mathbf{F}_q$.

**Lemma 9.1.** *The Teichmüller points of $A$ defined with respect to $\sigma$ are exactly those Teichmüller points of $R[X]^\dagger$ defined with respect to $\sigma$ which reduce modulo $\pi$ to the points of $X$.*

*Proof.* Let $\bar{x}$ be a geometric point of the variety $X$ of degree $d$ over $\mathbf{F}_q$. Let $R_d$ denote the unramified extension of $R$ with residue field $\mathbf{F}_{q^d}$. Recall that the Teichmüller lifting $x$ of $\bar{x}$ is the unique surjective $R$-algebra homomorphism $x \in \mathrm{Hom}_R(A, R_d)$ whose reduction mod $\pi$ is $\bar{x}$ such that $x \circ \sigma^d = x$. Composing $x$ with the reduction map $\mu$ (not the modulo $\pi$ reduction), we deduce that the composite map $x \circ \mu$ is a surjective $R$-algebra homomorphism in $\mathrm{Hom}_R(R[X]^\dagger, R_d)$ such that

$$\begin{aligned} (x \circ \mu) \circ \sigma^d &= x \circ (\mu \circ \sigma^d) \\ &= x \circ (\sigma^d \circ \mu) \\ &= (x \circ \sigma^d) \circ \mu \\ &= x \circ \mu \end{aligned}$$

where the first two $\sigma$ act on $R[X]^\dagger$ and the last two $\sigma$ act on $A$. By the uniqueness of the Teichmüller lifting, we deduce that $x \circ \mu$ is the Teichmüller lifting of the point on $R[X]^\dagger$ which reduces modulo $\pi$ to the point $\bar{x} \circ \bar{\mu}$ of $\mathbf{A}^n/\mathbf{F}_q$. The last point is just the point $\bar{x}$ of $X/\mathbf{F}_q$. By a cardinality counting of degree $d$ points, we deduce that the above procedure gives rise to all Teichmüller points of $R[X]^\dagger$ defined with respect to $\sigma$ which reduce modulo $\pi$ to the points of $X$. The lemma is proved. □

Next, we lift $(M, \phi)$ from $X/\mathbf{F}_q$ to $\mathbf{A}^n/\mathbf{F}_q$. Choose a row basis $\vec{e} = \{e_1, \cdots, e_r\}$ of $M$ over $A$ and write

$$M = \oplus_{i=1}^r A e_i, \quad \phi(\vec{e}) = \vec{e}C,$$

where $C$ is an $r \times r$ matrix with entries in $A$. Choose an $r \times r$ matrix $C_0$ with entries in $R[X]^\dagger$ such that $\mu(C_0) = C$. Then the relations

$$M_0 = \oplus_{i=1}^r R[X]^\dagger e_i, \quad \phi_0(\vec{e}) = \vec{e}C_0$$

define an overconvergent $\sigma$-module $(M_0, \phi_0)$ over $\mathbf{A}^n/\mathbf{F}_q$. The pair $(M_0, \phi_0)$ is our overconvergent lifting of $(M, \phi)$ from $X/\mathbf{F}_q$ to $\mathbf{A}^n/\mathbf{F}_q$. It is clear that the diagram

$$\begin{array}{ccccc} M_0 & \xrightarrow{\mu} & M & \longrightarrow & 0 \\ \downarrow \phi_0 & & \downarrow \phi & & \\ M_0 & \xrightarrow{\mu} & M & \longrightarrow & 0 \end{array}$$

commutes with exact rows, where $\mu$ acts coefficient-wise on $M_0$ and $\mu(e_i) = e_i$. By Lemma 9.1, we deduce the formula

(9.1) $$L(\phi/X, T) = L(\phi_0/X, T),$$

where the right side $L(\phi_0/X, T)$ is the L-function of the $\sigma$-module $\phi_0$ on $\mathbf{A}^n$ but the Euler product is now restricted to those Teichmüller points of $\mathbf{A}^n$ lying over points of $X$.

To reduce to the affine space case, we need to digress a little about Dwork's F-crystal, see [14] and [1] for further information. By going to a totally ramified finite extension of $R$ if necessary, we may assume that $R$ contains a primitive $p$-th root of unity. Let $\pi^*$ be an element in $R$ satisfying

$$(\pi^*)^{p-1} = -p.$$

Let

$$E(t) = \exp(\pi^*(t - t^q)).$$

By Dwork $p$-adic analytic construction of the additive character of $\mathbf{F}_q$, one knows that $E(t)$ is an overconvergent power series and $E(1)$ is a primitive $p$-th root of unity, see [14]. Let $\chi$ be the fixed non-trivial additive character of $\mathbf{F}_p$ such that $\chi(1) = E(1)$. For each positive integer $k$, let

$$\chi_k = \chi \circ \mathrm{Tr}_{\mathbf{F}_{q^k}/\mathbf{F}_p}.$$

This is a non-trivial additive character of the extension field $\mathbf{F}_{q^k}$. The order $p$ character $\chi$ defines a rank one overconvergent unit root F-crystal $E(t)$ over the affine $t$-line in Berthelot's sense, called Dwork's F-crystal.

Let $f_1(X), \cdots, f_s(X)$ be a set of generators for the ideal $I$. Its reduction modulo $\pi$ is then a system of defining equations for the variety $X$. Let

$$G_{s,f}(X, Y) \equiv Y_1 f_1(X) + \cdots + Y_s f_s(X) \pmod{\pi}$$

be the reduction modulo $\pi$, which is a polynomial over $\mathbf{F}_q$ in $n + s$ variables. Let

$$R[X, Y]^\dagger = R[X_1, \cdots, X_n, Y_1, \cdots, Y_s]^\dagger.$$

Extend $\sigma$ from $R[X]^\dagger$ to $R[X, Y]^\dagger$, for instance, by taking

$$Y_i^\sigma = \sigma(Y_i) = Y_i^q.$$

36 DAQING WAN

Let $\Phi_{s,f}(X,Y)$ be the $1 \times 1$ Frobenius matrix with respect to $\sigma$ of the pull-back of the overconvergent F-crystal $E(t)$ from the affine $t$-line $\mathbf{A}^1/\mathbf{F}_q$ to the affine space $\mathbf{A}^{n+s}/\mathbf{F}_q$ under the map

$$t \longrightarrow G_{s,f}(X,Y).$$

This gives a rank one unit root overconvergent F-crystal over $\mathbf{A}^{n+s}/\mathbf{F}_q$. Explicitly, we have

$$\Phi_{s,f}(X,Y) = \exp(\pi^* \sum_{i=1}^s Y_i f_i(X) - \pi^* \sum_{i=1}^s Y_i^\sigma f_i(X^\sigma)).$$

In particular, the power series $\Phi_{s,f}(X,Y)$ defines a rank one overconvergent unit root $\sigma$-module over the affine space $\mathbf{A}^{n+s}/\mathbf{F}_q$. At a Teichmüller point $(x,y)$ of $\mathbf{A}^{n+s}/\mathbf{F}_q$ of degree $k$, the $\sigma$-module $\Phi_{s,f}(X,Y)$ has the property that

$$\Phi_{s,f}(x,y) \cdots \Phi_{s,f}(x^{\sigma^{k-1}}, y^{\sigma^{k-1}}) = \chi_k(G_{s,f}(\bar{x}, \bar{y})),$$

where $\bar{x}$ and $\bar{y}$ denotes the reduction of $x$ and $y$ modulo $\pi$.

Via the inclusion $R[X]^\dagger \subset R[X,Y]^\dagger$, the $\sigma$-module $\phi_0$ extends to a $\sigma$-module over $\mathbf{A}^{n+s}/\mathbf{F}_q$, still denoted by $\phi_0$. The standard property of Dwork's function $E(t)$ implies the following result.

**Lemma 9.2.** *We have the relation*

$$(9.2) \qquad L(\phi_0/X, T) = L(\phi_0 \otimes \Phi_{s,f}/\mathbf{A}^{n+s}, \frac{1}{q^s} T).$$

*Proof.* The proof is the same as the standard method of expressing the zeta function of an affine variety in terms of the L-function of the exponential sum of a single polynomial. For the reader's convenience, we include a proof here. Write

$$L(\phi_0 \otimes \Phi_{s,f}/\mathbf{A}^{n+s}, T) = \exp(\sum_{k=1}^\infty \frac{S_k}{k} T^k).$$

Taking the logarithm on both sides and using the Euler product definition of L-function for the left side, one sees that $S_k$ is the $p$-adic character sum

$$S_k = \sum_{(x,y) \in T_k^{n+s}} \phi_0(x) \cdots \phi_0(x^{\sigma^{k-1}}) \Phi_{s,f}(x,y) \cdots \Phi_{s,f}(x^{\sigma^{k-1}}, y^{\sigma^{k-1}}),$$

where $T_k^{n+s}$ denotes the finite set of the Teichmüller liftings defined with respect to $\sigma$ of the rational points of $\mathbf{A}^{n+s}$ over the $k$-th extension field $\mathbf{F}_{q^k}$. Similarly, let $T_k^n$ (resp. $T_k^s$) denote the finite set of the Teichmüller liftings of rational points of $\bar{x} \in \mathbf{A}^n$ (resp. $\bar{y} \in \mathbf{A}^s$) over the $k$-th extension field $\mathbf{F}_{q^k}$. We can rewrite $S_k$ as

$$(9.3) \qquad S_k = \sum_{x \in T_k^n} \phi_0(x) \cdots \phi_0(x^{\sigma^{k-1}}) W_k(x),$$

where

$$W_k(x) = \sum_{y \in T_k^s} \Phi_{s,f}(x,y) \cdots \Phi_{s,f}(x^{\sigma^{k-1}}, y^{\sigma^{k-1}}).$$

At a Teichmüller point $(x,y)$ of $\mathbf{A}^{n+s}/\mathbf{F}_q$ of degree $k$, the $\sigma$-module $\Phi_{s,f}$ has the property that

$$\Phi_{s,f}(x,y) \cdots \Phi_{s,f}(x^{\sigma^{k-1}}, y^{\sigma^{k-1}}) = \chi_k(y_1 f_1(x) + \cdots + y_s f_s(x) \mod \pi).$$

Thus,
$$W_k(x) = \sum_{y_j \in \mathbf{F}_{q^k}} \chi_k(y_1 \bar{f}_1(\bar{x}) + \cdots + y_s \bar{f}_s(\bar{x})).$$

A standard argument shows that
$$W_k(x) = \begin{cases} q^{sk}, & \text{if } \bar{x} \in X, \\ 0, & \text{Otherwise.} \end{cases}$$

Substituting this formula into (9.3), we deduce that
$$S_k = q^{sk} \sum_{\bar{x} \in X(\mathbf{F}_{q^k})} \phi_0(x) \cdots \phi_0(x^{\sigma^{k-1}}),$$

where $x$ denotes the Teichmüller lifting of the point $\bar{x} \in X$ defined with respect to $\sigma$. Comparing this with the Euler product definition of $L(\phi_0, T)$, we conclude that (9.2) holds. The lemma is proved. □

Putting together equations (9.1)-(9.2), we conclude that the L-function of a $\sigma$-module over $X/\mathbf{F}_q$ can be expressed in terms of the L-function of a $\sigma$-module over the affine $(n+s)$-space $\mathbf{A}^{n+s}/\mathbf{F}_q$. Namely,

$$(9.4) \qquad L(\phi/X, T) = L(\phi_0 \otimes \Phi_{s,f}/\mathbf{A}^{n+s}, \frac{1}{q^s}T).$$

We will need a slightly more general construction involving tensor products and $k$-th powers. Let $(N, \psi)$ be another overconvergent $\sigma$-module over $X/\mathbf{F}_q$. In a similar way, we can lift $(N, \psi)$ to an overconvergent $\sigma$-module $(N_0, \psi_0)$ over $\mathbf{A}^n/\mathbf{F}_q$ and hence over $\mathbf{A}^{n+s}/\mathbf{F}_q$. The tensor product $(M_0 \otimes N_0, \phi_0 \otimes \psi_0)$ is a lifting of $(M \otimes N, \phi \otimes \psi)$ to an overconvergent $\sigma$-module over $\mathbf{A}^n/\mathbf{F}_q$ and hence over $\mathbf{A}^{n+s}/\mathbf{F}_q$. Applying (9.4) to $\psi_0 \otimes \phi_0^k$, we obtain

**Lemma 9.3.** *For each positive integer $k$, we have the formula*

$$(9.5) \qquad L(\psi \otimes \phi^k, T) = L(\psi_0 \otimes \phi_0^k \otimes \Phi_{s,f}/\mathbf{A}^{n+s}, \frac{1}{q^s}T).$$

*Proof.* Using the decomposition formula in (8.2) with $k_1 = 1$ and $k_2 = k$, the formula in (9.4) and then the decomposition formula in (8.2) again, we deduce that

$$\begin{aligned} L(\psi \otimes \phi^k, T) &= \prod_{i \geq 1} L(\psi \otimes \mathrm{Sym}^{k-i}\phi \otimes \wedge^i \phi, T)^{(-1)^{i-1}i} \\ &= \prod_{i \geq 1} L(\psi_0 \otimes \mathrm{Sym}^{k-i}\phi_0 \otimes \wedge^i \phi_0 \otimes \Phi_{s,f}/\mathbf{A}^{n+s}, \frac{1}{q^s}T)^{(-1)^{i-1}i} \\ &= L(\psi_0 \otimes \phi_0^k \otimes \Phi_{s,f}/\mathbf{A}^{n+s}, \frac{1}{q^s}T). \end{aligned}$$

The proof is complete. □

All the constructions and formulas of this section work perfectly well if we replace $A$ by $A_0$ and replace $R[X]^\dagger$ by the $\pi$-adic completion of $R[X]^\dagger$. The formula in (9.5) does not yet reduce our key lemma from $X$ to $\mathbf{A}^{n+s}$. The problem is that our lifted $\phi_0$ over $\mathbf{A}^n$ (and hence over $\mathbf{A}^{n+s}$) may not be ordinary at the slope zero side even if the initial $\phi$ over $X$ is ordinary at the slope zero side. To get around the difficulty, we shall apply the "pulling out" trick already used in our limiting approach.

**Lemma 9.4.** *Let $g$ be an invertible element of $A$. Then, there is a 1-unit $g_1$ in $A$ such that*
$$g^{q-1} = \frac{g^\sigma}{g} g_1^{q-1}, \quad g_1 = 1 + \pi g_2,$$
*where $g_2 \in A$.*

*Proof.* This is a version of Fermat's little theorem for the ring $A$. By our definition of the Frobenius map $\sigma$, we can write
$$g^q = g^\sigma + \pi h,$$
for some $h \in A$. Since $g$ is invertible in $A$, the element $g^\sigma$ is invertible in $A$ as well. Dividing the above equation by $g$, we obtain
$$g^{q-1} = \frac{g^\sigma}{g}(1 + \pi \frac{h}{g^\sigma}).$$
By the binomial theorem, one checks that
$$g_1 = (1 + \pi \frac{h}{g^\sigma})^{1/(q-1)}$$
is still a 1-unit in $A$. The proof is complete. □

Let now $(M, \phi)$ be an overconvergent $\sigma$-module over $X/\mathbf{F}_q$ ordinary at the slope zero side. Assume that the slope zero part $\phi(u)$ of $\phi$ has rank one. Let $\vec{e}$ be a row basis of $M$ over $A$ which is ordinary at the slope zero side. The matrix of $\phi$ under $\vec{e}$ has the form

(9.6) $$C = \begin{pmatrix} C_{00} & \pi C_{01} \\ C_{10} & \pi C_{11} \end{pmatrix},$$

where $C_{00}$ is an invertible element of $A$, and the other $C_{ij}$ are matrices over $A$.

**Definition 9.5.** We shall say that $C$ is in **normalized form** if $C_{00}$ is a 1-unit in $A$ and $C_{10}$ is divisible by $\pi$.

It is easy to see that if $C$ is in normalized form, then we can lift the overconvergent $\phi$ from $X/\mathbf{F}_q$ to an overconvergent $\phi_0$ over $\mathbf{A}^n/\mathbf{F}_q$ so that $\phi_0$ is still in normalized form and thus ordinary over $\mathbf{A}^n/\mathbf{F}_q$ fibre by fibre at the slope zero side. But it is not always possible to change a basis so that the matrix of $\phi$ is in normalized form. We shall proceed in two steps to reduce to a normalized situation.

First, we can change the basis so that $C_{10}$ is divisible by $\pi$. For this purpose, we consider the new basis
$$(e_1, e_2, \cdots, e_r) \begin{pmatrix} I_{00} & 0 \\ C_{10} C_{00}^{-1} & I_{11} \end{pmatrix} = \vec{e} E$$
of $M$ over $A$, where $I_{00}$ is the rank 1 identity matrix and $I_{11}$ is the $(r-1) \times (r-1)$ identity matrix. The transition matrix $E$ is overconvergent since $C_{10} C_{00}^{-1}$ is overconvergent. One calculates that the matrix of $\phi$ under the new basis $\vec{e}E$ is given by
$$E^{-1} C E^\sigma = \begin{pmatrix} C_{00} + \pi C_{01}(C_{10} C_{00}^{-1})^\sigma & \pi C_{01} \\ \pi (C_{11} - C_{10} C_{00}^{-1} C_{01})(C_{10} C_{00}^{-1})^\sigma & \pi(C_{11} - C_{10} C_{00}^{-1} C_{01}) \end{pmatrix}.$$
This matrix is still overconvergent since both $E$ and $E^\sigma$ are overconvergent. An important point is that the top left entry of this matrix is still a unit in $A$ and all other entries are now divisible by $\pi$.

Next, we show how to change the top left entry of the above matrix to be a 1-unit. Pulling out the invertible rank one factor $C_{00}$ from the matrix $E^{-1}CE^\sigma$, we see that $\phi$ can be written in the form

$$\phi = \xi \otimes \eta,$$

where $\eta$ is overconvergent in normalized form and $\xi$ is an overconvergent rank one unit root $\sigma$-module over $X$ with matrix $C_{00}$.

Let $\phi(u)$ be the slope zero part of $\phi$ and let $\eta(u)$ be the slope zero part of $\eta$. It is clear that we have the relation

$$\phi(u) = \xi \otimes \eta(u).$$

Thus, for each positive integer $k$, we get

$$\phi(u)^k = \xi^k \otimes \eta(u)^k,$$

where the $k$-th power is the same as the $k$-th tensor power since the rank is one. We shall need to consider the family $\phi(u)^k$ of rank one unit root $\sigma$-modules over $X$ parametrized by positive integer $k$. Let $m$ be the smallest non-negative residue of $k$ modulo $(q-1)$. That is,

$$k = m + (q-1)k_1,$$

for some integer $k_1$. We shall fix $m$ (there are only $q$ of them) and let $k_1$ vary. Since $C_{00}$ is invertible in $A$, Lemma 9.4 shows that there is a 1-unit $g_1$ in $A$ such that

$$C_{00}^{q-1} = \frac{\sigma(C_{00})}{C_{00}} g_1^{q-1}.$$

Thus,

$$\begin{aligned}
C_{00}^k &= C_{00}^m C_{00}^{(q-1)k_1} \\
&= C_{00}^m \frac{\sigma(C_{00}^{k_1})}{C_{00}^{k_1}} g_1^{(q-1)k_1} \\
&= \frac{C_{00}^m}{g_1^m} \frac{\sigma(C_{00}^{k_1})}{C_{00}^{k_1}} g_1^k.
\end{aligned}$$

Under a new basis, the matrix of $\xi^k$ becomes

$$\frac{C_{00}^m}{g_1^m} g_1^k.$$

Thus, we can write

$$\xi^k = \xi_1(m) \otimes \xi_2^k,$$

where $\xi_1(m)$ is an overconvergent rank one unit root $\sigma$-module over $X$ with matrix $C_{00}^m/g_1^m$ which depends only on $m$ not on $k_1$, and $\xi_2$ is a rank one overconvergent $\sigma$-module over $X$ whose matrix $g_1$ is a 1-unit in $A$. Thus,

$$\phi(u)^k = \xi_1(m) \otimes (\xi_2 \otimes \eta(u))^k = \xi_1(m) \otimes \rho^k,$$

where $\rho$ is the slope zero part of $\xi_2 \otimes \eta$. Note that $\xi_2 \otimes \eta$ is overconvergent in normalized form and ordinary at the slope zero side over $X$. Thus, we have shown that for all $k \equiv m \pmod{(q-1)}$,

(9.7) $$L(\psi \otimes \phi(u)^k, T) = L(\psi \otimes \xi_1(m) \otimes \rho^k, T),$$

where $\rho$ is the slope zero part of the overconvergent $\sigma$-module $\xi_2 \otimes \eta$ over $X/\mathbf{F}_q$ in normalized form and ordinary at the slope zero side. The key lemma is to prove the meromorphic continuation of the left side of (9.7). Replacing $\psi$ by $\psi \otimes \xi_1(m)$ and $\phi(u)$ by $\rho$, we deduce from (9.7) that our key lemma is reduced to the case that $\phi$ itself is already in normalized form. With this normalized assumption on $\phi$, we can now lift $(M, \phi)$ from $X$ to an overconvergent $(M_0, \phi_0)$ over $\mathbf{A}^n/\mathbf{F}_q$ (and hence over $\mathbf{A}^{n+s}/\mathbf{F}_q$) ordinary at the slope zero side and still in normalized form. The slope zero part $\phi_0(u)$ of $\phi_0$ is then a lifting of $\phi(u)$ to $\mathbf{A}^n/\mathbf{F}_q$ (and hence to $\mathbf{A}^{n+s}/\mathbf{F}_q$), not necessarily overconvergent.

Let $\psi_0$ be any overconvergent lifting of the overconvergent $\sigma$-module $\psi$ from $A$ to $R[X]^\dagger$, hence to $R[X,Y]^\dagger$. Applying Lemma 9.3 to the $\sigma$-module $\psi_0 \otimes \phi_0(u)^k$ over $\mathbf{A}^{n+s}/\mathbf{F}_q$, we obtain

$$(9.8) \qquad L(\psi \otimes \phi(u)^k, T) = L(\Phi_{s,f} \otimes \psi_0 \otimes \phi_0(u)^k/\mathbf{A}^{n+s}, \frac{1}{q^s}T).$$

The $\sigma$-module $\Phi_{s,f} \otimes \psi_0$ is clearly overconvergent over $\mathbf{A}^{n+s}/\mathbf{F}_q$. The rank one unit root $\sigma$-module $\phi_0(u)$ is the slope zero part of the overconvergent $\sigma$-module $\phi_0$ over $\mathbf{A}^{n+s}/\mathbf{F}_q$ ordinary at the slope zero side and in normalized form. It follows that our key lemma has been reduced from $X$ to the affine space case $\mathbf{A}^{n+s}$. More precisely, replacing $n+s$ by $n$ in the right side of (9.8), we have reduced our key lemma to the following result.

**Lemma 9.6.** *Let $\psi$ be an overconvergent $\sigma$-module over $\mathbf{A}^n/\mathbf{F}_q$. Let $\rho$ be a rank one unit root $\sigma$-module which is the slope zero part of some overconvergent $\sigma$-module $(M, \phi)$ over $\mathbf{A}^n/\mathbf{F}_q$ ordinary fibre by fibre at the slope zero side and in normalized form. Then the L-function $L(\psi \otimes \rho^k, T)$ is meromorphic everywhere for every integer $k$. Furthermore, the family of L-functions $L(\psi \otimes \rho^k, T)$ parametrized by integer $k$ is a strong family in the sense of* [23].

This lemma has already been proved in [19] in the case when the Frobenius lifting $\sigma$ is the simplest one $\sigma(X_i) = X_i^q$, where Dwork's trace formula can be applied. In particular, the more restrictive F-crystal version of Lemma 9.6 was already proved since the overconvergence condition is independent of the choice of $\sigma$. But for a general $\sigma$-module and a more general $\sigma$, we will need to use a sufficiently explicit form of Monsky's trace formula. The full Lemma 9.6, which is a special case of Theorem 7.7 in [23], is proved in [23].

## 10. Appendix: Proof Of Monsky Trace Formula

In this appendix, we give a proof of the extended Monsky trace formula for the finite rank case, based on various results in [13]. For this purpose, we need to define certain Grothendieck modules.

Recall that $A$ is our $n$-dimensional overconvergent integral domain and $\sigma$ is a fixed Frobenius endomorphism of $A$. Let $C(A, \sigma)$ be the category of finite $A$-modules $(M, \phi)$ with a $\sigma$-linear endomorphism $\phi$. Let $m(A, \sigma)$ be the free abelian group generated by the isomorphism classes of objects of $C(A, \sigma)$. Let $n(A, \sigma)$ be the subgroup of $m(A, \sigma)$ generated by the following two types of elements. The first type is of the form

$$(M, \phi) - (M_1, \phi_1) - (M_2, \phi_2),$$

where

$$0 \longrightarrow (M_1, \phi_1) \longrightarrow (M, \phi) \longrightarrow (M_2, \phi_2) \longrightarrow 0$$

is an exact sequence of objects in $C(A, \sigma)$. The second type is of the form

$$(M, \phi_1 + \phi_2) - (M, \phi_1) - (M, \phi_2),$$

where each of the three terms has the same underlying module $M$ and $\phi_1 + \phi_2$ simply denotes the sum of maps.

**Definition 10.1.** The Grothendieck group $K(A, \sigma)$ associated to the category $C(A, \sigma)$ is the quotient abelian group $m(A, \sigma)/n(A, \sigma)$.

If we started with the sub-category of finite projective $A$-modules (or even finite free $A$-modules) with a $\sigma$-linear map, its associated Grothendieck group is the same, by the existence of a finite projective resolution of a finite $A$-module. Thus, we can work with projective $\sigma$-modules if we wish. For $(M, \phi) \in C(A, \sigma)$, the class of $(M, \phi)$ in $K(A, \sigma)$ will be denoted by $[M, \phi]$. The abelian group $K(A, \sigma)$ has a natural $A$-module structure given by

$$a[M, \phi] = [M, L_a \circ \phi], \ a \in A,$$

where $L_a$ denotes the $A$-linear multiplication map by $a$. An object $(M, \phi)$ in $C(A, \sigma)$ is called trivial if $\phi$ is the zero map. It is clear that a trivial object $(M, \phi)$ in $C(A, \sigma)$ represents the zero element in $K(A, \sigma)$.

If $(M, \phi)$ is an object of $C(A, \sigma)$, the following diagram commutes:

$$\begin{array}{ccc} M & \xrightarrow{\phi \circ L_a} & M \\ \downarrow L_a & & \downarrow L_a \\ M & \xrightarrow{L_a \circ \phi} & M \end{array}$$

That is, the multiplication map $L_a$ induces a morphism

$$(M, \phi \circ L_a) \longrightarrow (M, L_a \circ \phi).$$

The kernel and cokernel of this morphism are trivial objects in $C(A, \sigma)$. Thus, they represent the zero element in $K(A, \sigma)$. This shows that

$$a[M, \phi] = [M, L_a \circ \phi] = [M, \phi \circ L_a] = [M, L_{\sigma(a)} \circ \phi] = \sigma(a)[M, \phi].$$

In particular,

$$(a - \sigma(a))[M, \phi] = 0$$

for every $a \in A$. Let $J$ be the ideal of $A$ generated by all the elements of the form $a - \sigma(a)$ with $a \in A$. Then, the quotient $A/J$ is isomorphic as $R$-algebra as a direct sum of $N(\bar{A})$ copies of $R$, where $N(\bar{A})$ denotes the number of $\mathbf{F}_q$-rational points on the affine $\mathbf{F}_q$-variety defined by $\bar{A}$. The above argument shows that the ideal $J$ annihilates $K(A, \sigma)$. Thus, the $A$-module $K(A, \sigma)$ has a natural $A/J$-module structure.

Similarly, let $C^*(A, \sigma)$ be the category of finite $A$-modules $(M, \Theta)$ with a Dwork operator $\Theta$:

$$\Theta(\sigma(a)m) = a\Theta(m), \ \ a \in A, m \in M.$$

The Grothendieck group $K^*(A, \sigma)$ associated to the category $C^*(A, \sigma)$ is defined in a similar way as $K(A, \sigma)$. The class of a Dwork operator $(M, \Theta)$ in $K^*(A, \sigma)$ will be denoted by $[M, \Theta]$. The abelian group $K^*(A, \sigma)$ has a natural $A$-module structure given by

$$a[M, \Theta] = [M, \Theta \circ L_a], \ a \in A.$$

Note that the definition of $A$-action on $K^*(A, \sigma)$ is different from that of $A$-action on $K(A, \sigma)$. An object $(M, \Theta)$ in $C^*(A, \sigma)$ is called trivial if $\Theta$ is the zero map. It is clear that a trivial object $(M, \Theta)$ in $C^*(A, \sigma)$ represents the zero element in $K^*(A, \sigma)$.

By a similar argument, one shows that the ideal $J$ annihilates $K^*(A, \sigma)$. Hence, the $A$-module $K^*(A, \sigma)$ has a natural $A/J$-module structure. Note that our notations are somewhat different from those in [13]. In particular, our $K^*(A, \sigma)$ (resp. $K(A, \sigma)$) is the $K(A, F)$ (resp. $K^*(A, F)$) in [13].

A basic result of Monsky (Theorem 3.5 and p.331 in [13]) is

**Lemma 10.2.** *The two modules $K(A, \sigma)$ and $K^*(A, \sigma)$ are both isomorphic to $A/J$ as $A/J$-modules.*

One can define several functorial maps from $K(A, \sigma)$ to $K^*(A, \sigma)$. Let $0 \leq i \leq n$ be an integer. Let $(M, \phi)$ be an object of $C(A, \sigma)$. Recall from section 3 that
$$M_i^* = \mathrm{Hom}_A(\Omega^i M, \Omega^n A) = \Omega^{n-i} M^\vee,$$
where
$$M^\vee = \mathrm{Hom}_A(M, A), \ \Omega^i M^\vee = M^\vee \otimes_A \Omega^i A.$$
A Dwork operator $\phi_i^*$ on $M_i^*$ is given by
$$\phi_i^*(f)(m) = \theta_n(f(\phi_i(m))), \ m \in \Omega^i M, \ f \in M_i^*,$$
where
$$\phi_i = \phi \otimes \sigma_i : \Omega^i M \longrightarrow \Omega^i M.$$
The corresponding Dwork operator on $\Omega^{n-i} M^\vee$ is denoted by $\theta_{n-i}(\phi)$. That is,
$$(M_i^*, \phi_i^*) = (\Omega^{n-i} M^\vee, \theta_{n-i}).$$
For each integer $0 \leq i \leq n$, we define a map $\lambda_i$ from the set $C(A, \sigma)$ to $K^*(A, \sigma)$ by
$$\lambda_i(M, \phi) = [M_i^*, \phi_i^*] = [\Omega^{n-i} M^\vee, \theta_{n-i}(\phi)].$$
The image $\lambda_i(M, \phi)$ depends only on the isomorphism class of the object $(M, \phi)$ in $C(A, \sigma)$. Thus, the map $\lambda_i$ extends to a homomorphism of abelian groups:
$$\lambda_i : m(A, \sigma) \longrightarrow K(A, \sigma).$$
It is easy to check that $\lambda_i$ annihilates the subgroup $n(A, \sigma)$ of $m(A, \sigma)$. It follows that $\lambda_i$ induces a well defined map from $K(A, \sigma)$ to $K^*(A, F)$, denoted by the same notation $\lambda_i$. This map respects the $A$-module (resp. $A/J$-module) structure:
$$\lambda_i(a[M, \phi]) = a\lambda_i([M, \phi]), \ a \in A$$
by our definition of the $A$-module structure of $K(A, \sigma)$ and $K^*(A, \sigma)$. We deduce

**Lemma 10.3.** *For each integer $0 \leq i \leq n$, the map $\lambda_i : K(A, \sigma) \to K^*(A, \sigma)$ is a homomorphism of cyclic $A/J$-modules.*

The identity element $[A, \sigma]$ is a generator of the cyclic $A/J$-module $K(A, \sigma)$, denoted by $\chi(A)$. It image under $\lambda_0$ is
$$\lambda_0(\chi(A)) = \lambda_0([A, \sigma]) = [A_0^*, \sigma_0^*] = [\Omega^n A, \theta_n].$$
This is the identity element of the cyclic $A/J$-module $K^*(A, \sigma)$, denoted by $\chi^*(A)$. Thus, Lemma 10.3 can be strengthened as follows (see p.331 in [13]).

**Lemma 10.4.** *The map $\lambda_0 : K(A, \sigma) \to K^*(A, \sigma)$ is an isomorphism of cyclic $A/J$-modules.*

For $i > 0$, the map $\lambda_i$ will not be an isomorphism but it can be shown to be topologically nilpotent. Define a new map $\lambda$ from $K(A,\sigma)$ to $K^*(A,F)$ by

$$\lambda = \sum_{i=0}^{n}(-1)^i \lambda_{n-i}.$$

That is,

$$\begin{aligned}
\lambda([M,\phi]) &= \sum_{i=0}^{n}(-1)^i \lambda_{n-i}([M,\phi]) \\
&= \sum_{i=0}^{n}(-1)^i [M^*_{n-i}, \phi^*_{n-i}] \\
&= \sum_{i=0}^{n}(-1)^i [\Omega^i M^\vee, \theta_i(\phi)].
\end{aligned}$$

The map $\lambda$ is a homomorphism of cyclic $A/J$-modules. We have the following stronger lemma:

**Lemma 10.5.** *The map $\lambda : K(A,\sigma) \to K^*(A,\sigma)$ is an isomorphism of cyclic $A/J$-modules.*

This result is a consequence of the next lemma which shows that $\lambda$, like $\lambda_0$, also takes the generator $\chi(A)$ of $K(A,\sigma)$ to the generator $\chi^*(A)$ of $K^*(A,\sigma)$.

**Lemma 10.6.** *As elements of $K^*(A,\sigma)$,*

$$\lambda(\chi(A)) = \chi^*(A).$$

*Proof.* In our notations, Theorem 5.1 in [13] says that

$$\chi^*(A) = \sum_{i=0}^{n}(-1)^i [\Omega^i A, \theta_i].$$

Since the dual $\sigma$-module of the identity element $(A,\sigma)$ is itself,

$$\begin{aligned}
\lambda(\chi(A)) &= \lambda([A,\sigma]) \\
&= \sum_{i=0}^{n}(-1)^i [A^*_{n-i}, \sigma^*_{n-i}] \\
&= \sum_{i=0}^{n}(-1)^i [\Omega^i A, \theta_i] \\
&= \chi^*(A).
\end{aligned}$$

The lemma is true. $\square$

Lemma 10.6 is essentially the rank one case of Monsky's trace formula. To extend it to higher rank case, we need to introduce a trace map on the Grothendieck module $K(A,\sigma)$. If $(M,\phi)$ is an object of $C(A,\sigma)$, then the quotient $M/JM$ is a finite $A/J$-module. The induced action of $\phi$ on $M/JM$ is $A/J$-linear. The trace of the $A/J$-linear map acting on $M/JM$ is called the $A/J$-trace of $\phi$, denoted by $\text{Tr}_{A/J}(M,\phi)$ or in short by $\text{Tr}_{A/J}(\phi)$. It depends only on the isomorphism class of $(M,\phi)$ in $C(A,\sigma)$ and hence extends to an $A$-linear map on $m(A,\sigma)$. It is easy to see that the $A/J$-trace map annihilates the subgroup $n(A, \sigma)$ and thus induces a well defined $A/J$-linear map on the Grothendieck module $K(A, \sigma)$:

$$\text{Tr}_{A/J}: \; K(A, J) \longrightarrow A/J,$$

denoted by the same notation. Since the $A/J$-trace of the identity object $(A, \sigma)$ is 1, which is a generator of $A/J$, we obtain

**Lemma 10.7.** *The $A/J$-trace map $\text{Tr}_{A/J}: K(A, J) \to A/J$ is an isomorphism of cyclic $A/J$-modules.*

Since $K(A, \sigma)$ is a free $A/J$-module of rank 1 with generator $\chi(A)$, for a given object $(M, \phi) \in C(A, \sigma)$, there is a unique element $a(\phi) \in A/J$ such that

$$[M, \phi] = a(\phi)\chi(A).$$

Taking the $A/J$-trace on both sides, we obtain the relation

$$\text{Tr}_{A/J}(\phi) = a(\phi).$$

In other words, we have

**Lemma 10.8.** *As elements of $K(A, \sigma)$, we have*

$$[M, \phi] = \text{Tr}_{A/J}(\phi)\chi(A).$$

Applying the $\lambda$ map to both sides of the above equation and using Lemma 10.6, we obtain the main result of this section:

**Theorem 10.9.** *Let $(M, \phi) \in C(A, \sigma)$. As elements of $K^*(A, \sigma)$, we have*

$$\begin{aligned}
\text{Tr}_{A/J}(\phi)\chi^*(A) &= \sum_{i=0}^{n}(-1)^i[M^*_{n-i}, \phi^*_{n-i}] \\
&= \sum_{i=0}^{n}(-1)^i[\Omega^i M^\vee, \theta_i(\phi)].
\end{aligned}$$

This is the abstract form of the general trace formula. To obtain the concrete trace formula, we let Tr denote the $R$-trace map from $K^*(A, \sigma)$ to $R$. Taking the $R$-trace in the equation of Theorem 10.9, we obtain

**Theorem 10.10.** *Let $(M, \phi) \in C(A, \sigma)$. Then,*

$$\begin{aligned}
\sum_{x} x(\text{Tr}_{A/J}(\phi)) &= \sum_{i=0}^{n}(-1)^i \text{Tr}(\phi^*_{n-i}|M^*_{n-i}) \\
&= \sum_{i=0}^{n}(-1)^i \text{Tr}(\theta_i(\phi)|\Omega^i M^\vee),
\end{aligned}$$

*where $x$ runs over the set of Teichmüller points of degree 1, that is the set of surjective $R$-algebra homomorphisms $x: A \to R$ such that $x \circ \sigma = x$.*

This theorem is the additive form of Theorem 3.1 in section 3. It is true for free $M$ by our discussion above. The general case follows by taking a free resolution of $(M, \phi)$.

# References

[1] P. Berthelot, *Géométrie rigide et cohomologie des variétés algébriques de caractéristique p*, Mém. Soc. Math. France (N.S.) **23**(1986), 3, 7–32.

[2] P. Berthelot, *Finitude et pureté cohomologique en cohomologie rigide*, Invent. Math., **128** (1997), no. 2, 329–377.

[3] P. Deligne, *La Conjecture de Weil*, II, Publ. Math. IHES, **52**(1980), 137–252.

[4] A.J. de Jong, *Barsotti-Tate groups and crystals*, Proceedings ICM 1998(Berlin), Vol. 2, 259–265.

[5] B. Dwork, *Normalized period matrices*, Ann. Math., **94**(1971), 337–388.

[6] B. Dwork, *Normalized period matrices* II, Ann. Math., **98**(1973), 1–57.

[7] J.-Y. Etesse and B. Le Stum, *Fonctions L associées aux F-isocristaux surconvergents* I, Interprétation cohomologique, Math. Ann., **296**(1993), 557–576.

[8] W. Fulton, *A note on weakly complete algebras*, Bull. Amer. Math. Soc., **75**(1969), 591–593.

[9] N. Katz, *Travaux de Dwork*, Séminaire Bourbaki, exposé 409(1971/72), Lecture Notes in Math., No. 317, 1973, 167–200.

[10] N. Katz, *Slope filtration of F-crystals*, Astérisque, **63**(1979), 113–164.

[11] B. Mazur, *Frobenius and the Hodge filtration*, Bull. Amer. Math. Soc., **78**(1972), 653–667.

[12] Z. Mebkhout, *Sur le thorme de finitude de la cohomologie p-adique d'une varit affine non singulire*, Amer. J. Math., **119**(1997), 1027–1088.

[13] P. Monsky, *Formal cohomology* III, Ann. Math., **93**(1971), 315–343.

[14] P. Mosky, *p-Adic Analysis and Zeta Functions*, Kinokukuniya Book-Store, 1970.

[15] M. van der Put, *The cohomology of Monsky and Washnitzer*, Mém. Soc. Math. France (N.S.) **23**(1986), 33–60.

[16] J-P, Serre, *Endomorphismes complétements continues des espaces de Banach p-adiques*, Publ. Math., IHES., **12**(1962), 69–85.

[17] D. Wan, *Noetherian subrings of power series rings*, Proc. Amer. Math. Soc., **123**(1995), 1681–1686.

[18] D. Wan, *Meromorphic continuation of L-functions of p-adic representations*, Ann. Math., **143**(1996), 469–498.

[19] D. Wan, *Dwork's conjecture on unit root zeta functions*, Ann. Math., **150**(1999), 867–927.

[20] D. Wan, *A quick introduction to Dwork's conjecture*, Contemporary Mathematics, **245**(1999), 147–163.

[21] D. Wan, *Pure L-functions from algebraic geometry over finite fields*, Proc. of the Augsburg Conference on Finite Fields, to appear.

[22] D. Wan, *Poles of zeta functions of complete intersections*, Chinese Ann. Math., to appear.

[23] D. Wan, *Rank one case of Dwork's conjecture*, J. Amer. Math. Soc., this issue.